\documentclass[thmsa,onecolumn,10pt]{article}
\usepackage{amssymb}
\usepackage{amsfonts}
\usepackage{graphicx}
\usepackage{amsmath}

\newtheorem{theorem}{Theorem}

\newtheorem{lemma}[theorem]{Lemma}

\newtheorem{proposition}[theorem]{Proposition}

\newenvironment{proof}[1][Proof]{\textbf{#1.} }{\ \rule{0.5em}{0.5em}}

\begin{document}

\title{Martingale Structure of Skorohod Integral Processes}
\author{Giovanni PECCATI \\
Laboratoire de Statistique Th\'{e}orique et Appliqu\'{e}e\\
Universit\'{e} Paris VI, France\\
email: giovanni.peccati@gmail.com\\\\
Mich\`{e}le THIEULLEN \\
Laboratoire de Probabilit\'{e}s et Mod\`{e}les Al\'{e}atoires\\
Universit\'{e}s Paris VI and Paris VII, France\\
email: mth@ccr.jussieu.fr\\\\
Ciprian A. TUDOR\footnote{The new address of C.A. Tudor is:
SAMOS/MATISSE Universit\'{e} de Pantheon - Sorbonne, Paris I, France.}\\
Laboratoire de Probabilit\'{e}s et Mod\`{e}les Al\'{e}atoires\\
Universit\'{e}s Paris VI and Paris VII, France\\
email: tudor@ccr.jussieu.fr\\}
\date{February 3, 2005}
\maketitle

\begin{abstract}
Let the process $\left\{ Y_{t}, t\in[0,1] \right\}$, have the form $Y_{t}=\delta \left( u%
\mathbf{1}_{\left[ 0,t\right] }\right) $, where $\delta $ stands for a
Skorohod integral with respect to Brownian motion, and $u$ is a measurable
process verifying some suitable regularity conditions. We use a recent
result by Tudor (2004), to prove that $Y_{t}$ can be represented as the
limit of linear combinations of processes that are products of forward and
backward Brownian martingales. Such a result is a further step towards the
connection between the theory of continuous-time (semi)martingales, and that
of anticipating stochastic integration. We establish an explicit link
between our results and the classic characterization, due to Duc and Nualart
(1990), of the chaotic decomposition of Skorohod integral processes. We also
explore the case of Skorohod integral processes that are time-reversed
Brownian martingales, and provide an \textquotedblleft
anticipating\textquotedblright\ counterpart to the classic Optional Sampling
Theorem for It\^{o} stochastic integrals.



\textbf{Running title --}Martingale structure of integrals
\end{abstract}

\section{Introduction}

Let $\left( C_{\left[ 0,1\right] },\mathcal{C},\mathbb{P}\right) =\left(
\Omega ,\mathcal{F},\mathbb{P}\right) $ be the canonical space, where $%
\mathbb{P}$ is the law of a standard Brownian motion started from zero, and
write $X=\left\{ X_{t}:t\in \left[ 0,1\right] \right\} $ for the coordinate
process. In this paper, we investigate some properties of \textit{Skorohod
integral processes} defined with respect to $X$, that is, measurable
stochastic processes with the form
\begin{equation}
Y_{t}=\int_{0}^{1}u_{s}\mathbf{1}_{\left[ 0,t\right] }\left(
s\right) dX_{s}=\int_{0}^{t}u_{s}dX_{s},\ \ \ t\in \left[
0,1\right] ,  \label{SkoInt}
\end{equation}%
where $\left\{ u_{s}:s\in \left[ 0,1\right] \right\} $ is a suitably regular
(and not necessarily adapted) process verifying
\begin{equation}
\mathbb{E}\left[ \int_{0}^{1}u_{s}^{2}ds\right] <+\infty \text{,}
\label{integrability1}
\end{equation}%
and the stochastic differential $dX$ has to be interpreted in the
Skorohod sense (as defined in Skorohod (1975); see the discussion
below, as well as Nualart and Pardoux (1988) or Nualart (1995,
Chapters 1 and 3), for basic results concerning Skorohod
integration). It is well known that if $u_{s}$ is adapted to the
natural filtration of $X$ (noted $\left\{ \mathcal{F}_{s}:s\in
\left[ 0,1\right] \right\} $) and satisfies
(\ref{integrability1}), then $Y_{t}$ is a stochastic integral
process in the It\^{o} sense (as defined e.g. in Revuz and Yor
(1999)), and therefore $Y_{t}$ is a square-integrable
$\mathcal{F}_{t}$ - martingale. In
general, the martingale property of $Y_{t}$ fails when $u_{s}$ is not $%
\mathcal{F}_{s}$ - adapted, and $Y_{t}$ may have a path behavior
that is very different from the ones of classical It\^{o}
stochastic integrals (see Barlow and Imkeller (1992), for examples
of anticipating integral processes with very irregular
trajectories). However, in Tudor (2004) it is proved that the
class of Skorohod integral processes (when the integrand $u$ is
sufficiently regular) coincides with the set of
\textit{Skorohod-It\^{o} integrals}, i.e. processes admitting the
representation
\begin{equation}
Y_{t}=\int_{0}^{t}\mathbb{E}\left[ v_{s}\mid \mathcal{F}_{\left[ s,t\right]
^{c}}\right] dX_{s}\text{, \ \ }t\in \left[ 0,1\right] \text{,}
\label{ISproc}
\end{equation}%
where $v$ is measurable and satisfies (\ref{integrability1}), $\mathcal{F}_{%
\left[ s,t\right] ^{c}}:=\mathcal{F}_{s}\vee \sigma \left\{
X_{1}-X_{r}:r\geq t\right\} $, and for each fixed $t$ the stochastic
integral is in the usual It\^{o} sense (indeed, for fixed $t$, $X_{s}$ is a
standard Brownian motion on $\left[ 0,t\right] $, with respect to the
enlarged filtration $s\mapsto \mathcal{F}_{\left[ s,t\right] ^{c}}$).

The principal aim of this paper is to use representation (\ref{ISproc}), in
order to provide an exhaustive characterization of Skorohod integral
processes in terms of products of \textit{forward} and \textit{backward}
Brownian martingales. In particular, we shall prove that a process $Y_{t}$
has the representation (\ref{SkoInt}) (or, equivalently, (\ref{ISproc})) if,
and only if, $Y_{t}$ is the limit, in an appropriate norm, of linear
combinations of stochastic processes of the type
\begin{equation*}
Z_{t}=M_{t}\times N_{t}\text{, \ \ }t\in \left[ 0,1\right] \text{,}
\end{equation*}%
where $M_{t}$ is a centered (forward) $\mathcal{F}_{t}$ - martingale, and $%
N_{t}$ is a $\mathcal{F}_{\left[ 0,t\right] ^{c}}$ - backward martingale
(that is, for any $0\leq s<t\leq 1$, $N_{t}\in \mathcal{F}_{\left[ 0,t\right]
^{c}}$ and $\mathbb{E}\left[ N_{s}\mid \mathcal{F}_{\left[ 0,t\right] ^{c}}%
\right] =N_{t}$ ). Such a representation accounts in particular
for the well-known property of Skorohod integral processes (see
e.g. Nualart (1995, Lemma 3.2.1):
\begin{equation}
\mathbb{E}\left[ Y_{t}-Y_{s}\mid \mathcal{F}_{\left[ s,t\right] ^{c}}\right]
=0\mbox{ for every }s<t,  \label{MartSkoh}
\end{equation}%
playing in the anticipating calculus a somewhat analogous role as
the martingale property in the It\^{o}'s calculus. We will see, in
the subsequent discussion, that our characterization of processes
such as $Y_{t}$ complements some classic results contained in Duc
and Nualart (1990), where the authors study the multiple Wiener
integral expansion of Skorohod integral processes.

The paper is organized as follows. In Section 2, we introduce some
notation and discuss preliminary issues concerning the Malliavin
calculus; in Section 3, the main results of the paper are stated
and proved; in Section 4, we establish an explicit link between
our results and those contained in Duc and Nualart (1990); in
Section 5, we concentrate on a special class of Skorohod
integral processes, whose elements can be represented as \textit{%
time-reversed Brownian martingales}, and we state sufficient conditions to
have that such processes are semimartingales in their own filtration;
eventually, Section 6 discusses some relations between processes such as (%
\ref{SkoInt}) and stopping times.

\section{Notation and preliminaries}

Let $L^{2}\left( \left[ 0,1\right] ,dx\right) =L^{2}\left( \left[ 0,1\right]
\right) $ be the Hilbert space of square integrable functions on $\left[ 0,1%
\right] .$ In what follows, the notation
\begin{equation*}
X=\left\{ X\left( f\right) :f\in L^{2}\left( \left[ 0,1\right] \right)
\right\}
\end{equation*}%
will indicate an \textit{isonormal Gaussian process }on \thinspace $%
L^{2}\left( \left[ 0,1\right] \right) $, that is, $X$ is a centered Gaussian
family indexed by the elements of $L^{2}\left( \left[ 0,1\right] \right) $,
defined on some (complete) probability space $\left( \Omega ,\mathcal{F},%
\mathbb{P}\right) $ and such that $\mathbb{E}\left[ X\left( f\right) X\left(
g\right) \right] =\int_{0}^{1}f\left( x\right) g\left( x\right) dx$ for
every $f,g\in L^{2}\left( \left[ 0,1\right] \right) $. We also introduce the
standard Brownian motion $X_{t}=X\left( \mathbf{1}_{\left[ 0,t\right]
}\right) $, $t\in \left[ 0,1\right] $, and note $L^{2}\left( \mathbb{P}%
\right) $ the space of square integrable functionals of $X$. The
usual notation of Malliavin calculus is adopted throughout the
sequel (see Nualart (1990)): for instance, $D$ and $\delta $
denote the (Malliavin) derivative operator and the Skorohod
integral with respect to the Wiener process $X$. For $k\geq 1$ and
$p\geq 2$, $\mathbb{D}^{k,p}$ denotes the
space of $k$ times differentiable functionals of $X$, endowed with the norm $%
\left\Vert \cdot \right\Vert _{k,p}$, whereas $\mathbb{L}^{k,p}=L^{p}\left( %
\left[ 0,1\right] ;\mathbb{D}^{k,p}\right) $. Note that $\mathbb{L}%
^{k,p}\subset Dom(\delta )$, the domain of $\delta $. Now take a
Borel subset $A$ of $[0,1]$, and denote by $\mathcal{F}_{A}$ the
$\sigma $-field generated by random variables with the form
$X\left( f\right) $, where $f\in L^{2}\left( \left[ 0,1\right]
\right) $ is such that its support is contained in $A$. We recall
that if $F\in {\mathcal{F}}_{A}$ and $F\in \mathbb{D}^{1,2}$, then
\begin{equation}
D_{t}F\left( \omega \right) =0\mbox{ on }A^{c}\times \Omega .  \label{d}
\end{equation}%
We will also need the following integration by parts formula:
\begin{equation}
\delta (Fu)=F\delta (u)-\int_{[0,1]}D_{s}Fu_{s}ds  \label{ip}
\end{equation}%
p.s. - $\mathbb{P}$, whenever $u\in Dom(\delta )$ and $F\in \mathbb{D}^{1,2}$
are such that $\mathbb{E}(F^{2}\int_{[0,1]}u_{s}^{2}ds)<\infty $.

\medskip

Eventually, let us introduce, for further reference, the following families
of $\sigma $-fields:
\begin{eqnarray*}
\mathcal{F}_{t} &=&\sigma \left\{ X_{h}:h\leq t\right\} \text{, \ \ }t\in %
\left[ 0,1\right] ; \\
\mathcal{F}_{\left[ s,t\right] ^{c}} &=&\sigma \left\{ X_{h}:h\leq s\right\}
\vee \sigma \left\{ X_{1}-X_{h}:h\geq t\right\} \text{, \ \ }0\leq s<t\leq 1,
\end{eqnarray*}%
and observe that, to simplify the notation, we will write $\mathcal{F}_{%
\left[ 0,t\right] ^{c}}=\mathcal{F}_{t^{c}}$, so that $\mathcal{F}_{\left[
s,t\right] ^{c}}=\mathcal{F}_{t^{c}}\vee \mathcal{F}_{s}$.

\section{Skorohod integral processes and martingales}

Let $L_{0}^{2}\left( \mathbb{P}\right) $ denote the space of zero mean
square integrable functionals of $X$. We write $Y\in \mathbf{BF}$ to
indicate that the measurable stochastic process $Y=\left\{ Y_{t}:t\in \left[
0,1\right] \right\} $ can be represented as a finite linear combination of
processes with the form%
\begin{equation}
Z_{t}=\mathbb{E}\left[ H_{1}\mid \mathcal{F}_{t}\right] \times \mathbb{E}%
\left[ H_{2}\mid \mathcal{F}_{t^{c}}\right] =M_{t}\times N_{t},\text{ \ \ }%
t\in \left[ 0,1\right] \text{,}  \label{FB}
\end{equation}%
where $H_{1}\in L_{0}^{2}\left( \mathbb{P}\right) $ and $H_{2}\in
L^{2}\left( \mathbb{P}\right) $. Note that $M$ in (\ref{FB}) is a forward
(centered) Brownian martingale, whereas $N$ is a backward Brownian
martingale. For every measurable process $G=$ $\left\{ G_{t}:t\in \left[ 0,1%
\right] \right\} $, we also introduce the notation%
\begin{equation}
V\left( G\right) =\sup_{\pi }\mathbb{E}\left[ \sum_{j=0}^{m-1}\left(
G_{t_{j}}-G_{t_{j+1}}\right) ^{2}\right] \text{,}  \label{Varia}
\end{equation}%
where $\pi $ runs over all partitions of $\left[ 0,1\right] $ with the form $%
0=t_{0}<t_{1}<...<t_{m}=1$. The following result shows that $\mathbf{BF}$ is
in some sense dense in the class of Skorohod integral processes.

\begin{theorem}
Let $u\in \mathbb{L}^{k,p}$, with $k\geq 3$ and $p>2$. Then, there exists a
sequence of processes
\begin{equation*}
\left\{ Z_{t}^{\left( r\right) }:t\in \left[ 0,1\right] \right\} ,\text{ \ \
}r\geq 1,
\end{equation*}%
with the following properties:

\begin{description}
\item[(i)] for every $r$, $Z^{\left( r\right) }\in \mathbf{BF}$;

\item[(ii)] for every $r$, $Z_{t}^{\left( r\right) }=\int_{0}^{t}\mathbb{E}%
\left[ v_{\alpha }^{\left( r\right) }\mid \mathcal{F}_{\left[ \alpha ,t%
\right] ^{c}}\right] dX_{\alpha }$, $t\in \left[ 0,1\right] $, where $%
v^{\left( r\right) }\in \mathbb{L}^{k-2,p}$;

\item[(iii)] for every $r$, $V\left( Z^{\left( r\right) }\right) <+\infty $
and $\lim_{r\rightarrow \infty }V\left( \delta \left( u\mathbf{1}_{\left[
0,\cdot \right] }\right) -Z^{\left( r\right) }\right) =0$.
\end{description}
\end{theorem}

\bigskip

Note that points (i) and (iii) of Theorem 1 imply that $Z^{\left(
r\right) }$ converges to $\delta \left( u\mathbf{1}_{\left[
0,\cdot \right] }\right) $ uniformly in $L^{2}\left(
\mathbb{P}\right) $. This implies that the convergence takes also
place in the sense of finite dimensional distributions. Before
proving Theorem 1, we need to state two simple results.

\begin{lemma}
Fix $k\geq 1$ and $p\geq 2$. Let $A_{1}$ and $A_{2}$ be two disjoint subsets
of $\left[ 0,1\right] $, and let $\mathcal{F}_{A_{i}}$, $i=1,2$, be the $%
\sigma $-field generated by random variables of the form $X(h\mathbf{1}_{A_{i}})$, $%
h\in L^{2}\left(\left[ 0,1\right] \right) $. Suppose that $F\in \mathcal{F}_{A_{1}}\vee \mathcal{F}%
_{A_{2}} $ and also $F\in \mathbb{D}^{k,p}$. Then, $F$ is the limit in $%
\mathbb{D}^{k,p}$ of linear combinations of smooth random variables of the
type%
\begin{equation}
G=G_{1}\times G_{2}\text{,}  \label{product}
\end{equation}%
where, for $i=1,2$, $G_{i}$ is smooth and $\mathcal{F}_{A_{i}}$ - measurable.
\end{lemma}

\begin{proof}
By definition, every $F\in \mathbb{D}^{k,p}$ can be approximated in the
space $\mathbb{D}^{k,p}$ by a sequence of smooth polynomial functionals of
the type%
\begin{equation*}
P_{m}=p_{n_{\left( m\right) }}\left( X\left( h_{1}^{\left( m\right) }\right)
,...,X\left( h_{n_{\left( m\right) }}^{\left( m\right) }\right) \right)
\text{, \ \ }m\geq 1\text{,}
\end{equation*}%
where, for every $m$, $n_{\left( m\right) }\geq 1$, $p_{n_{\left( m\right)
}} $ is a polynomial in $n_{\left( m\right) }$ variables and, for $%
j=1,...,n_{\left( m\right) }$, $h_{j}^{\left( m\right) }\in L^{2}\left( %
\left[ 0,1\right] \right) $. It is also easily checked that $\mathbb{E}\left[
P_{m}\mid \mathcal{F}_{A_{1}}\vee \mathcal{F}_{A_{2}}\right] \in \mathbb{D}%
^{k,p}$ for every $m$ and, since $F\in \mathcal{F}_{A_{1}}\vee \mathcal{F}%
_{A_{2}}$,
\begin{equation*}
\mathbb{E}\left[ P_{m}\mid \mathcal{F}_{A_{1}}\vee \mathcal{F}_{A_{2}}\right]
\rightarrow F
\end{equation*}%
in $\mathbb{D}^{k,p}$. To conclude, it is sufficient to prove that every
random variable of the kind
\begin{equation*}
Z=\mathbb{E}\left[ \left( X\left( h_{1}\right) \right) ^{k_{1}}\cdot \cdot
\cdot \left( X\left( h_{n}\right) \right) ^{k_{n}}\mid \mathcal{F}%
_{A_{1}}\vee \mathcal{F}_{A_{2}}\right]
\end{equation*}%
where $h_{j}\in L^{2}\left( \left[ 0,1\right] \right) $ and $k_{j}\geq 1$,
can be represented as a linear combination of random variables such as (\ref%
{product}). To see this, write $A_{3}=\left[ 0,1\right] \backslash \left(
A_{1}\cup A_{2}\right) $, and use twice the binomial formula to obtain%
\begin{eqnarray*}
\left( X\left( h_{j}\right) \right) ^{k_{j}} &=&\sum_{l=0}^{k_{j}}\binom{%
k_{j}}{l}\left( X\left( h_{j}\mathbf{1}_{A_{1}}\right) \right)
^{k_{j}-l}\left( X\left( h_{j}\mathbf{1}_{A_{2}\cup A_{3}}\right) \right)
^{l} \\
&=&\sum_{l=0}^{k_{j}}\sum_{a=0}^{l}\binom{k_{j}}{l}\binom{l}{a}\left(
X\left( h_{j}\mathbf{1}_{A_{1}}\right) \right) ^{k_{j}-l}\left( X\left( h_{j}%
\mathbf{1}_{A_{2}}\right) \right) ^{l-a}\left( X\left( h_{j}\mathbf{1}%
_{A_{3}}\right) \right) ^{a}\text{,}
\end{eqnarray*}%
thus implying that the functional $\left( X\left( h_{1}\right) \right)
^{k_{1}}\cdot \cdot \cdot \left( X\left( h_{n}\right) \right) ^{k_{n}}$ is a
linear combination of random variables of the type
\begin{equation*}
H=\prod_{j=1}^{n}\left( X\left( h_{j}\mathbf{1}_{A_{1}}\right) \right)
^{\gamma _{1,j}}\left( X\left( h_{j}\mathbf{1}_{A_{2}}\right) \right)
^{\gamma _{2,j}}\left( X\left( h_{j}\mathbf{1}_{A_{3}}\right) \right)
^{\gamma _{3,j}}\text{,}
\end{equation*}%
where $\gamma _{i,j}\geq 0$, $j=1,...,n$, $i=1,2,3$. To conclude, use
independence to obtain
\begin{equation*}
\mathbb{E}\left[ H\mid \mathcal{F}_{A_{1}}\vee \mathcal{F}_{A_{2}}\right] =%
\mathbb{E}\left[ \prod_{j=1}^{n}\left( X\left( h_{j}\mathbf{1}%
_{A_{3}}\right) \right) ^{\gamma _{3,j}}\right] \times \prod_{j=1}^{n}\left(
X\left( h_{j}\mathbf{1}_{A_{1}}\right) \right) ^{\gamma
_{1,j}}\prod_{j=1}^{n}\left( X\left( h_{j}\mathbf{1}_{A_{2}}\right) \right)
^{\gamma _{2,j}},
\end{equation*}%
and therefore the desired conclusion.
\end{proof}

\bigskip

\textbf{Remark -- }Suppose that $F=I_{n}^{X}\left( h\right) $, $n\geq 1$,
where $I_{n}^{X}$ stands for a multiple Wiener integral of order $n$. Then, $%
F\in \mathbb{D}^{k,p}$ for every $k\geq 1$ and $p\geq 2$. Moreover, the
isometric properties of multiple integrals imply that $F$ can be
approximated in $\mathbb{D}^{k,2}$, and therefore in $\mathbb{D}^{k,p}$ for
every $p\geq 2$, by linear combinations of random variables with the form $%
H_{n}\left( X\left( h\right) \right) $, where $H_{n}$ is an Hermite
polynomial of the $n$th order and $h$ is an element of $L^{2}\left( \left[
0,1\right]\right) $. In particular, if $F\in \mathcal{F}_{A_{1}}\vee
\mathcal{F}_{A_{2}}$ as in the statement of Lemma 2, the arguments contained
in the above proof entail that $F$ is the limit in $\mathbb{D}^{k,p}$ of
linear combinations of random variables of the type $G=G_{1}\times G_{2}$,
where, for $i=1,2$, $G_{i}$ is a $\mathcal{F}_{A_{i}}$ - measurable
polynomial functional of order $\gamma _{i}\geq 0$ such that $\gamma
_{1}+\gamma _{2}\leq n$.

\bigskip

The proof of the following result is trivial, and it is therefore omitted.

\bigskip

\begin{lemma}
Fix $k\geq 1$ and $p\geq 2$, as well as a partition $%
0=t_{0}<t_{1}<...<t_{n}=1$ of $\left[ 0,1\right] $.\ Then, for every finite
collection $\left\{ F_{j}:j=1,...,n\right\} $ of elements of $\mathbb{D}%
^{k,p}$, the process
\begin{equation*}
u_{t}=\sum_{j=0}^{n-1}F_{j}\mathbf{1}_{\left( t_{j},t_{j+1}\right) }\left(
t\right)
\end{equation*}%
is an element $\mathbb{L}^{k,p}.$ Moreover, if $F_{j}^{m}\underset{%
m\rightarrow +\infty }{\rightarrow }F_{j}$ in $\mathbb{D}^{k,p}$, then, as $%
m\rightarrow +\infty $, the sequence of processes%
\begin{equation*}
u_{t}^{m}=\sum_{j=0}^{n-1}F_{j}^{m}\mathbf{1}_{\left( t_{j},t_{j+1}\right)
}\left( t\right)
\end{equation*}%
converges to $u$ in $\mathbb{L}^{k,p}$.
\end{lemma}

\bigskip

\noindent \textbf{Proof of Theorem 1. }It is well known (see e.g. Duc and
Nualart (1990)) that the process $t\mapsto Y_{t}=\delta \left( u%
\mathbf{1}_{\left[ 0,t\right] }\right) $ is such that $V\left(
Y\right) <+\infty $. Moreover, according to Proposition 1 in Tudor
(2004), $Y$ admits the (unique) representation
\begin{equation}
Y_{t}=\int_{0}^{t}\mathbb{E}\left[ v_{\alpha }\mid \mathcal{F}_{\left[
\alpha ,t\right] ^{c}}\right] dX_{\alpha },\text{ \ \ }t\in \left[ 0,1\right]
,  \label{TudorRepres}
\end{equation}%
where $v\in \mathbb{L}^{k-2,p}$. Now, for every partition $\pi $ of the type
$0=t_{0}<...<t_{n}=1$, we introduce the step process
\begin{equation}
v_{t}^{\pi }=\sum_{i=0}^{n-1}\frac{1}{t_{i+1}-t_{i}}\left(
\int_{t_{i}}^{t_{i+1}}\mathbb{E}\left[ v_{s}\mid \mathcal{F}_{\left[
t_{i},t_{i+1}\right] ^{c}}\right] ds\right) \mathbf{1}_{\left(
t_{i},t_{i+1}\right) }\left( t\right) ,\text{ \ }t\in \left[ 0,1\right]
\text{,}  \label{CondStrato}
\end{equation}%
and we recall that $v^{\pi }\in \mathbb{L}^{k-2,p}$, and that $v^{\pi }$
converges to $v$ in $\mathbb{L}^{k-2,p}$ whenever the mesh of $\pi $, noted $%
\left\vert \pi \right\vert $, converges to zero. Now define
$Y_{t}^{\pi }=\int_{0}^{t}\mathbb{E}\left[ v_{\alpha }^{\pi }\mid
\mathcal{F}_{\left[ \alpha ,t\right] ^{c}}\right] dX_{\alpha }$.
From the calculations contained in Tudor (2004, Proposition 2), we
deduce that
\begin{equation}
V\left( Y-Y^{\pi }\right) \leq \left\Vert v-v^{\pi }\right\Vert _{1,2}^{2},
\label{Inequa}
\end{equation}%
and therefore that $V\left( Y^{\pi }\right) <+\infty $ and $V\left( Y-Y^{\pi
}\right) $ converges to zero, as $\left\vert \pi \right\vert \rightarrow 0$.
Now fix a partition $\pi $, and note, for $i=0,...,n-1$,
\begin{equation}
F_{i}^{\pi }:=\frac{1}{t_{i+1}-t_{i}}\left( \int_{t_{i}}^{t_{i+1}}\mathbb{E}%
\left[ v_{s}\mid \mathcal{F}_{\left[ t_{i},t_{i+1}\right] ^{c}}\right]
ds\right) \in \mathcal{F}_{\left[ t_{i},t_{i+1}\right] ^{c}}.
\label{RVapproaching}
\end{equation}%
Since for every $i$ and every $s$ such that $t_{i}\leq s\leq t_{i+1}$ and $%
s<t$,%
\begin{equation*}
\mathbb{E}\left[ F_{i}^{\pi }\mid \mathcal{F}_{\left[ s,t\right] ^{c}}\right]
=\mathbb{E}\left[ F_{i}^{\pi }\mid \mathcal{F}_{\left[ s,t\right] ^{c}\cap %
\left[ t_{i},t_{i+1}\right] ^{c}}\right] =\mathbb{E}\left[ F_{i}^{\pi }\mid
\mathcal{F}_{\left[ t_{i},t_{i+1}\vee t\right] ^{c}}\right] \text{ ,}
\end{equation*}%
we obtain, using the properties (\ref{ip}) and (\ref{d})%
\begin{eqnarray*}
Y_{t}^{\pi } &=&\sum_{i=0}^{n-1}\int_{0}^{t}\mathbf{1}_{\left[ t_{i},t_{i+1}%
\right] }\left( s\right) \mathbb{E}\left[ F_{i}^{\pi }\mid \mathcal{F}_{%
\left[ t_{i},t_{i+1}\vee t\right] ^{c}}\right] dX_{s} \\
&=&\sum_{i=0}^{n-1}\mathbb{E}\left[ F_{i}^{\pi }\mid \mathcal{F}_{\left[
t_{i},t_{i+1}\vee t\right] ^{c}}\right] \left( X_{t\wedge
t_{i+1}}-X_{t_{i}}\right) \mathbf{1}_{\left( t\geq t_{i}\right) } \\
&=&\sum_{i=0}^{n-1}Z_{t}^{\left( \pi ,i\right) }\text{,}
\end{eqnarray*}%
where $Z_{t}^{\left( \pi ,i\right) }=\mathbb{E}\left[ F_{i}^{\pi }\mid
\mathcal{F}_{\left[ t_{i},t_{i+1}\vee t\right] ^{c}}\right] \left(
X_{t\wedge t_{i+1}}-X_{t_{i}}\right) \mathbf{1}_{\left( t\geq t_{i}\right) }$%
. Now fix $i=0,...,n-1$. Since $F_{i}^{\pi }$ is $\mathcal{F}_{\left[
t_{i},t_{i+1}\right] ^{c}}$-measurable and $F_{i}\in \mathbb{D}^{k-2,p}$,
thanks to Lemma 2 in the special case $A_{1}=\left( 0,t_{i}\right) $ and $%
A_{2}=\left( t_{i+1},1\right) $, the random variable $F_{i}^{\pi }$ is the
limit in the space $\mathbb{D}^{k-2,p}$ of a sequence of random variables of
the type%
\begin{equation}
G_{m}^{\left( i,\pi \right) }=\sum_{k=1}^{M_{m}}G_{m,k}^{\left( i,\pi
,1\right) }\times G_{m,k}^{\left( i,\pi ,2\right) }\text{, \ \ }m\geq 1\text{%
,}  \label{approxi}
\end{equation}%
where, for every $m$, $M_{m}\geq 1$ and, for every $k$, $G_{m,k}^{\left(
i,\pi ,1\right) },G_{m,k}^{\left( i,\pi ,2\right) }$ are smooth and such
that $G_{m,k}^{\left( i,\pi ,1\right) }\in \mathcal{F}_{t_{i}}$, and $%
G_{m,k}^{\left( i,\pi ,2\right) }\in \mathcal{F}_{t_{i+1}^{c}}$. This
implies, thanks to Lemma 3, that the process
\begin{equation*}
v_{t}^{m,\pi }=\sum_{i=0}^{n-1}G_{m}^{\left( i,\pi \right) }\mathbf{1}%
_{\left( t_{i},t_{i+1}\right) }\left( t\right) ,\text{ \ }t\in \left[ 0,1%
\right] \text{,}
\end{equation*}%
converges to $v^{\pi }$ in $\mathbb{L}^{k-2,p}$, and therefore, due to an
inequality similar to (\ref{Inequa}), for every $\pi $ the sequence of
processes%
\begin{eqnarray*}
Y_{t}^{m,\pi } &=&\sum_{i=0}^{n-1}\int_{0}^{t}\mathbf{1}_{\left[
t_{i},t_{i+1}\right] }\left( s\right) \mathbb{E}\left[ G_{m}^{\left( i,\pi
\right) }\mid \mathcal{F}_{\left[ t_{i},t_{i+1}\vee t\right] ^{c}}\right]
dX_{s} \\
&=&\sum_{i=0}^{n-1}\mathbb{E}\left[ G_{m}^{\left( i,\pi \right) }\mid
\mathcal{F}_{\left[ t_{i},t_{i+1}\vee t\right] ^{c}}\right] \left(
X_{t\wedge t_{i+1}}-X_{t_{i}}\right) \mathbf{1}_{\left( t\geq t_{i}\right) }
\\
&=&\sum_{i=0}^{n-1}\sum_{k=1}^{M_{m}}\mathbb{E}\left[ G_{m,k}^{\left( i,\pi
,1\right) }\times G_{m,k}^{\left( i,\pi ,2\right) }\mid \mathcal{F}_{\left[
t_{i},t_{i+1}\vee t\right] ^{c}}\right] \left( X_{t\wedge
t_{i+1}}-X_{t_{i}}\right) \mathbf{1}_{\left( t\geq t_{i}\right) } \\
&=&\sum_{i=0}^{n-1}\sum_{k=1}^{M_{m}}U_{t}^{\left( m,k,\pi ,i\right) }\text{%
, \ \ }m\geq 1\text{,}
\end{eqnarray*}%
is such that $V\left( Y^{m,\pi }\right) <+\infty $ and $\lim_{m\rightarrow
+\infty }V\left( Y^{\pi }-Y^{m,\pi }\right) =0$. We shall now show that $%
U^{\left( m,k,\pi ,i\right) }\in \mathbf{BF}$. As a matter of fact,
\begin{eqnarray}
U_{t}^{\left( m,k,\pi ,i\right) } &=&\mathbb{E}\left[ G_{m,k}^{\left( i,\pi
,1\right) }G_{m,k}^{\left( i,\pi ,2\right) }\mid \mathcal{F}_{\left[
t_{i},t_{i+1}\vee t\right] ^{c}}\right] \left( X_{t\wedge
t_{i+1}}-X_{t_{i}}\right) \mathbf{1}_{\left( t\geq t_{i}\right) }
\label{FmartBmart} \\
&=&\left[ G_{m,k}^{\left( i,\pi ,1\right) }\left( X_{t\wedge
t_{i+1}}-X_{t_{i}}\right) \mathbf{1}_{\left( t\geq t_{i}\right) }\right]
\times \mathbb{E}\left[ G_{m,k}^{\left( i,\pi ,2\right) }\mid \mathcal{F}_{%
\left[ t_{i},t_{i+1}\vee t\right] ^{c}}\right]   \notag \\
&=&M_{t}\times N_{t}.  \notag
\end{eqnarray}%
Eventually, observe that $M_{t}=\int_{0}^{t}H_{s}dX_{s}$ where $%
H_{s}=G_{m,k}^{\left( i,\pi ,1\right) }\mathbf{1}_{\left(
t_{i},t_{i+1}\right) }\left( s\right) $, and therefore, since $H_{s}$ is $%
\mathcal{F}_{s}$ - predictable, $M_{t}$ is a Brownian martingale such that $%
M_{0}=0$; on the other hand,
\begin{eqnarray}
N_{t} &=&\mathbb{E}\left[ G_{m,k}^{\left( i,\pi ,2\right) }\mid \mathcal{F}_{%
\left[ t_{i},t_{i+1}\vee t\right] ^{c}}\right] =\mathbb{E}\left[ \mathbb{E}%
\left[ G_{m,k}^{\left( i,\pi ,2\right) }\mid \mathcal{F}_{t_{i+1}^{c}}\right]
\mid \mathcal{F}_{\left[ t_{i},t_{i+1}\vee t\right] ^{c}}\right]
\label{Bmart} \\
&=&\mathbb{E}\left[ \mathbb{E}\left[ G_{m,k}^{\left( i,\pi ,2\right) }\mid
\mathcal{F}_{t_{i+1}^{c}}\right] \mid \mathcal{F}_{\left( t_{i+1}\vee
t\right) ^{c}}\right] =\mathbb{E}\left[ G_{m,k}^{\left( i,\pi ,2\right)
}\mid \mathcal{F}_{\left( t_{i+1}\vee t\right) ^{c}}\right] \text{,}  \notag
\end{eqnarray}%
and also%
\begin{eqnarray}
N_{t} &=&\mathbb{E}\left[ \mathbb{E}\left[ G_{m,k}^{\left( i,\pi ,2\right)
}\mid \mathcal{F}_{t_{i+1}^{c}}\right] \mid \mathcal{F}_{\left( t_{i+1}\vee
t\right) ^{c}}\right]   \label{Bmart2} \\
&=&\mathbb{E}\left[ \mathbb{E}\left[ G_{m,k}^{\left( i,\pi ,2\right) }\mid
\mathcal{F}_{t_{i+1}^{c}}\right] \mid \mathcal{F}_{t^{c}}\right]   \notag \\
&=&\mathbb{E}\left[ N_{0}\mid \mathcal{F}_{t^{c}}\right] ,  \notag
\end{eqnarray}%
so that $N_{t}$ is a backward martingale such that
$N_{1}=\mathbb{E}\left[ G_{m,k}^{\left( i,\pi ,2\right) }\right]
$. As a consequence, we obtain that $U^{\left( m,k,\pi ,i\right)
}$, and therefore $Y^{m,\pi }$, is an element of $\mathbf{BF}$. We
have therefore shown that for every $r\geq 1$ there exists a
partition $\pi \left( r\right) $ and a number $m\left( r,\pi
\left( r\right) \right) $ such that $V\left( Y-Y^{\pi \left(
r\right) }\right) \leq 1/\left( 4r\right) $ and also $V\left(
Y^{\pi \left( r\right) ,m\left( r,\pi \left( r\right) \right)
}-Y^{\pi \left( r\right) }\right) \leq 1/\left( 4r\right) $. To
conclude, set $Z^{\left( r\right) }:=Y^{\pi \left( r\right)
,m\left( r,\pi \left( r\right) \right) }$ and observe that
\begin{equation*}
V\left( Y-Z^{\left( r\right) }\right) \leq 2\left[ V\left( Y-Y^{\pi \left(
r\right) }\right) +V\left( Y^{\pi \left( r\right) ,m\left( r,\pi \left(
r\right) \right) }-Y^{\pi \left( r\right) }\right) \right] \leq \frac{1}{r}%
\text{.}
\end{equation*}

$\blacksquare $

\bigskip

The next result contains a converse to Theorem 1.

\begin{theorem}
Let the sequence $Z^{\left( n\right) }\in \mathbf{BF}$, $n\geq 1$, be such
that $V\left( Z^{\left( n\right) }\right) <+\infty $ and
\begin{equation*}
\lim_{n,m\rightarrow +\infty }V\left( Z^{\left( n\right) }-Z^{\left(
m\right) }\right) =0.
\end{equation*}%
Then, there exists a process $\left\{ Y_{t}:t\in \left[ 0,1\right] \right\} $
such that

\begin{description}
\item[(i)] $Y_{t}$ admits a Skorohod integral representation;

\item[(ii)] $V\left( Y\right) <+\infty $ and $\lim_{n\rightarrow +\infty
}V\left( Z^{\left( n\right) }-Y\right) =0.$
\end{description}
\end{theorem}

\begin{proof}
We shall first prove point (ii). Consider the trivial partition $t_{0}=0$, $%
t_{1}=1$. Then, the assumptions in the statement (remember that $%
Z_{0}^{\left( n\right) }=0$) imply that $Z_{1}^{\left( n\right) }$ is a
Cauchy sequence in $L^{2}\left( \mathbb{P}\right) $. Moreover, since for
every $t\in \left( 0,1\right) $,
\begin{equation*}
\lim_{n,m\rightarrow +\infty }\mathbb{E}\left[ \left(
Z_{t}^{\left( n\right) }-Z_{t}^{\left( m\right) }\right)
^{2}+\left( Z_{t}^{\left( n\right) }-Z_{t}^{\left( m\right)
}-\left( Z_{1}^{\left( n\right) }-Z_{1}^{\left( m\right) }\right)
\right) ^{2}\right] =0,
\end{equation*}%
we readily obtain that for every $t\in \left[ 0,1\right] $ there exists $%
Y_{t}\in L^{2}\left( \mathbb{P}\right) $ such that $Y_{0}=0$ and also $%
Z_{t}^{\left( n\right) }\rightarrow Y_{t}$ in $L^{2}\left( \mathbb{P}\right)
$. Now fix $\varepsilon >0$; it follows from the assumptions that there
exists $N\geq 1$ such that for every $n,m>N$ and for every partition $%
0=t_{0}<...<t_{M}=1$%
\begin{equation*}
\mathbb{E}\left[\sum_{j=0}^{M-1}\left( \left( Z_{t_{j+1}}^{\left(
n\right) }-Z_{t_{j+1}}^{\left( m\right) }\right) -\left(
Z_{t_{j}}^{\left( n\right) }-Z_{t_{j}}^{\left( m\right) }\right)
\right) ^{2}\right]\leq \varepsilon ,
\end{equation*}%
and therefore, letting $m$ go to infinity, we obtain that for $n>N$
\begin{equation*}
\sup_{\pi }\mathbb{E}\left[\sum_{j=0}^{M-1}\left( \left(
Z_{t_{j+1}}^{\left( n\right) }-Y_{t_{j+1}}\right) -\left(
Z_{t_{j}}^{\left( n\right) }-Y_{t_{j}}\right) \right)
^{2}\right]=V\left( Z^{\left( n\right) }-Y\right) \leq \varepsilon
,
\end{equation*}%
that entails $\lim_{n\rightarrow +\infty }V\left( Z^{\left( n\right)
}-Y\right) =0$. To conclude the proof of (ii), observe that, for $n>N$ as
before,
\begin{equation*}
V\left( Y\right) \leq 2\left( V\left( Y-Z^{\left( n\right) }\right) +V\left(
Z^{\left( n\right) }\right) \right) \leq 2\left( \varepsilon +V\left(
Z^{\left( n\right) }\right) \right) <+\infty .
\end{equation*}

Thanks to Proposition 2.3. in Duc and Nualart (1990), to show
point
(i) it is now sufficient to prove that for any $s<t$%
\begin{equation*}
\mathbb{E}\left[ Y_{t}-Y_{s}\mid \mathcal{F}_{\left[ s,t\right] ^{c}}\right]
=0,
\end{equation*}%
which is easily proven by using $L^{2}$ convergence as well as the fact that
for every process $Z_{t}$ as in (\ref{FB}) we have%
\begin{equation*}
\mathbb{E}\left[ Z_{t}-Z_{s}\mid \mathcal{F}_{\left[ s,t\right] ^{c}}\right]
=N_{t}\mathbb{E}\left[ M_{t}\mid \mathcal{F}_{\left[ s,t\right] ^{c}}\right]
-M_{s}\mathbb{E}\left[ N_{s}\mid \mathcal{F}_{\left[ s,t\right] ^{c}}\right]
=0.
\end{equation*}
\end{proof}

\section{Representation of finite chaos Skorohod integral processes}

We say that the process $Y=\left\{ Y_{t}:t\in \left[ 0,1\right]
\right\} $ is a \textit{finite chaos Skorohod integral process of
order }$N\geq 0$ (written: $Y\in \mathbf{FS}_{N}$) if
$Y_{t}=\delta \left( u\mathbf{1}_{\left[ 0,t\right] }\right) $ for
some Skorohod integrable process $u_{\alpha }\left( \omega \right)
\in L^{2}\left( \left[ 0,1\right] \times \Omega \right) $ such
that, for each $\alpha \in \left[ 0,1\right] $, the random
variable $u_{\alpha }$ belongs to $\oplus _{j=0,...,N}C_{j}$,
where $C_{j}$ represents the $j$th Wiener chaos associated to $X$.
Note that if $Y\in \mathbf{FS}_{N}$, then, for each $t$, $Y_{t}\in
\oplus _{j=0,...,N+1}C_{j}$. We also define $\mathbf{FS}=\cup
_{N\geq 0}\mathbf{FS}_{N}$. The aim of this paragraph is to
discuss the relations between the results of the previous section,
and the representation of the elements of the class $\mathbf{FS}$
introduced in Duc and Nualart (1990)). To this end, we shall need
some further notation (note that our formalism is essentially
analogous to the one contained in the first part of Duc and
Nualart (1990)).

\bigskip

For every $M\geq 2$ and every $1\leq m\leq M$, we write $\mathbf{j}_{\left(
m\right) }\subset \{1,...,M\} $ to indicate that the vector $%
\mathbf{j}_{\left( m\right) }=\left( j_{1},...,j_{m}\right) $ has
integer-valued components such that $1\leq
j_{1}<j_{2}<...<j_{m}\leq M$.
Note that $\mathbf{j}_{\left( M\right) }=\left( 1,...,M\right) $. We set $%
\mathbf{j}_{\left( 0\right) }=\varnothing $ by definition, and also, given $%
\mathbf{x}_{M}=\left( x_{1},...,x_{M}\right) \in \left[ 0,1\right] ^{M}$ and
$\mathbf{j}_{\left( m\right) }=\left( j_{1},...,j_{m}\right) \subset \{1,...,M\} ,$%
\begin{equation*}
\mathbf{x}_{\mathbf{j}_{\left( m\right) }}:=\left(
x_{j_{1}},...,x_{j_{m}}\right) \text{ \ \ ; \ \ }\mathbf{x}_{\mathbf{j}%
_{\left( 0\right) }}:=0.
\end{equation*}%
We use the following notation: (a) for every permutation $\sigma
^{M}=\{\sigma \left( 1\right) ,...,\sigma \left( M\right) \}$ of
$\{1,...,M\},$ we set%
\begin{equation*}
\Delta _{M}^{\sigma ^{M}}:=\left\{ \left( x_{1},...,x_{M}\right) \in \left[
0,1\right] ^{M}:0<x_{\sigma \left( M\right) }<...<x_{\sigma \left( 1\right)
}<1\right\}
\end{equation*}%
and also write%
\begin{equation*}
\Delta _{M}^{\sigma _{0}^{M}}:=\Delta _{M}=\left\{ \left(
x_{1},...,x_{M}\right) \in \left[ 0,1\right] ^{M}:0<x_{M}<...<x_{1}<1\right\}
\end{equation*}%
for the simplex contained in $\left[ 0,1\right] ^{M}$; (b) for every $%
m=0,...,M$ and $\mathbf{j}_{\left( m\right) }\subset \{1,...,M\}$,%
\begin{equation*}
\Delta _{M}^{\mathbf{j}_{\left( m\right) }}:=\left\{ \left(
x_{1},...,x_{M}\right) \in \left( 0,1\right) ^{M}:\max_{i\in \mathbf{j}%
_{\left( m\right) }}\left( x_{i}\right) <\min_{l\in \{1,...,M\}
\backslash \mathbf{j}_{\left( m\right) }}\left( x_{l}\right)
\right\} ,
\end{equation*}%
where $\max_{i\in \mathbf{\varnothing }}\left( x_{i}\right) :=0$ and $%
\min_{l\in \varnothing }\left( x_{l}\right) :=1$; (c) for every $t\in \left[
0,1\right] $ and every $\mathbf{j}_{\left( m\right) }\subset \{1,...,M\},$%
\begin{equation*}
\Delta _{M}^{\mathbf{j}_{\left( m\right) }}\left( t\right) :=\left\{ \left(
x_{1},...,x_{M}\right) \in \left( 0,1\right) ^{M}:\max_{i\in \mathbf{j}%
_{\left( m\right) }}\left( x_{i}\right) <t<\min_{l\in \{1,...,M\}
\backslash \mathbf{j}_{\left( m\right) }}\left( x_{l}\right)
\right\} ;
\end{equation*}%
(d) for every $t\in \left[ 0,1\right] $,%
\begin{equation*}
A_{M,m}\left( t\right) =\bigcup\limits_{\mathbf{j}_{\left(
m\right) }\subset \{1,...,M\}}\Delta _{M}^{\mathbf{j}_{\left(
m\right) }}\left( t\right) .
\end{equation*}

\textbf{Remark -- }Note that $\Delta _{M}^{\mathbf{j}_{\left( 0\right) }}=$ $%
\Delta _{M}^{\mathbf{j}_{\left( M\right) }}=\left( 0,1\right)
^{M}$ and, in general, for every $m=0,...,M$ and every
$\mathbf{j}_{\left(
m\right) }\subset \{1,...,M\}$%
\begin{equation*}
\Delta _{M}^{\mathbf{j}_{\left( m\right) }}=\bigcup\limits_{t\in \mathbb{Q}%
\cap \left( 0,1\right) }\Delta _{M}^{\mathbf{j}_{\left( m\right) }}\left(
t\right) .
\end{equation*}

We have also the following relations,
\begin{equation*}
A_{M,M}\left( t\right) =\Delta _{M}^{\mathbf{j}_{\left( M\right) }}\left(
t\right) =\left( 0,t\right) ^{M}\text{ \ ; \ }A_{M,0}\left( t\right) =\Delta
_{M}^{\mathbf{j}_{\left( 0\right) }}\left( t\right) =\left( t,1\right) ^{M}%
\text{,}
\end{equation*}%
and moreover, if $t\in \left\{ 0,1\right\} $ and $0<m<M$, then $%
A_{M,m}\left( t\right) =\varnothing $.

\bigskip

The following result corresponds to properties $\left( \mathbf{B1}\right) $-$%
\left( \mathbf{B3}\right) $ in Duc and Nualart (1990).

\begin{proposition}
Fix $M\geq 2$ and $0\leq m\leq M$, and let the previous notation prevail.
Then, (i)
\begin{equation*}
\bigcup\limits_{\mathbf{j}_{\left( m\right) }\subset
\{1,...,M\}}\Delta _{M}^{\mathbf{j}_{\left( m\right) }}=\left[
0,1\right] ^{M},\text{ \ \ a.e.-}Leb\text{,}
\end{equation*}%
where $Leb$ stands for Lebesgue measure; (ii) if
$\mathbf{i}_{\left( m\right) },\mathbf{j}_{\left( m\right)
}\subset \{1,...,M\}$,
then $\Delta _{M}^{\mathbf{j}_{\left( m\right) }}\cap \Delta _{M}^{\mathbf{i}%
_{\left( m\right) }}\neq \varnothing $ if, and only if, $\mathbf{i}_{\left(
m\right) }=\mathbf{j}_{\left( m\right) }$; (iii) for any $t\in \left[ 0,1%
\right] $, if $m\neq m^{\prime }$ and $0\leq m,m^{\prime }\leq M$, then $%
A_{M,m}\left( t\right) \cap A_{M,m^{\prime }}\left( t\right) =\varnothing $,
and also%
\begin{equation*}
\bigcup\limits_{m=0,...,M}A_{M,m}\left( t\right) =\left[ 0,1\right] ^{M},%
\text{ \ \ a.e.-}Leb\text{.}
\end{equation*}
\end{proposition}

\bigskip

The next fact is a combination of Theorems 1.3 and 2.1 in Duc and
Nualart (1990), and gives a univocal characterization of the chaos
expansion of the elements of $\mathbf{FS}$. Note that, in the
following, we will write $L_{s}^{2}\left( \left[ 0,1\right]
^{k}\right) $, $k\geq 2$, to indicate the set of symmetric
functions on $\left[ 0,1\right] ^{k}$ that are square integrable
with respect to Lebesgue measure. Moreover,
for any $k\geq 2$ and $f\in L_{s}^{2}\left( \left[ 0,1\right] ^{k}\right) $, the symbol $%
I_{k}^{X}\left( f\right) $ will denote the standard multiple Wiener-It%
\^{o} integral (of order $k$) of $f$ with respect to $X$ (see e.g.
Nualart (1995, 1998) for definitions). We will also use the
notation $L_{s}^{2}\left( \left[ 0,1\right]
\right) =L^{2}\left( \left[ 0,1\right] \right) $ and, for $f\in L^{2}\left( %
\left[ 0,1\right] \right) $, $I_{1}^{X}\left( f\right) =X\left(
f\right) $.

\begin{theorem}[Duc and Nualart]
Let the above notation prevail, and fix $N\geq 0$. Then, the process $%
Y=\left\{ Y_{t}:t\in \left[ 0,1\right] \right\} $ is an element of $\mathbf{%
FS}_{N}$ if, and only if, there exists a (unique) collection of kernels $%
\left\{ f_{l,q}:1\leq q\leq l\leq N+1\right\} $ such that $f_{l,q}\in
L_{s}^{2}\left( \left[ 0,1\right] ^{l}\right) $ for every $1\leq q\leq l\leq
N+1$ and%
\begin{equation}
Y_{t}=\sum_{l=1}^{N+1}\sum_{q=1}^{l}I_{l}^{X}\left( f_{l,q}\mathbf{1}%
_{A_{l,q}\left( t\right) }\right) \text{, \ \ }t\in \left[ 0,1\right] .
\label{FSN}
\end{equation}%
Moreover, if condition (\ref{FSN}) is satisfied
\begin{equation}
\sum_{l=1}^{N+1}l!\sum_{q=0}^{l-1}\left\Vert f_{l,q}-f_{l,q+1}\right\Vert
^{2}\leq V\left( Y\right) <+\infty ,  \label{majoration}
\end{equation}%
where $V\left( Y\right) $ is defined according to (\ref{Varia}), and $%
f_{l,0}:=0$.
\end{theorem}

\bigskip

The link between the objects introduced in this paragraph and those of the
previous section is given by the following

\begin{lemma}
Fix $m,n\geq 0$, and for every $r\geq 1$ take a natural number $M_{r}\geq 1$%
, as well as two collections of kernels
\begin{equation*}
\left\{ h_{j}^{\left( u,r\right) }:1\leq u\leq M_{r};j=1,...,m\right\} \text{
\ ; \ }\left\{ g_{i}^{\left( u,r\right) }:1\leq u\leq
M_{r};i=1,...,n\right\} ,
\end{equation*}%
where $h_{j}^{\left( u,r\right) }\in L_{s}^{2}\left( \left[ 0,1\right]
^{j}\right) $ and $g_{i}^{\left( u,r\right) }\in L_{s}^{2}\left( \left[ 0,1%
\right] ^{i}\right) $ for every $i,j$, and a set of real numbers
\begin{equation*}
\left\{ b^{\left( u,r\right) }:1\leq u\leq M_{r}\right\} .
\end{equation*}%
For every $t\in \left[ 0,1\right] $ and $r\geq 1$, we define%
\begin{equation}
Z_{t}^{\left( r\right) }:=\sum_{u=1}^{M_{r}}Z_{t}^{\left( u,r\right)
}=\sum_{u=1}^{M_{r}}\left( \sum_{j=1}^{m}I_{j}^{X}\left( h_{j}^{\left(
u,r\right) }\mathbf{1}_{\left( 0,t\right) }^{\otimes j}\right) \right)
\times \left( b^{\left( u,r\right) }+\sum_{i=1}^{n}I_{i}^{X}\left(
g_{i}^{\left( u,r\right) }\mathbf{1}_{\left( t,1\right) }^{\otimes i}\right)
\right) .  \label{approaching}
\end{equation}%
Then: (i) for every $r\geq 1$, $V\left( Z^{\left( r\right) }\right) <+\infty
$; (ii) if
\begin{equation*}
\lim_{r,r^{\prime }\uparrow +\infty }V\left( Z^{\left( r\right) }-Z^{\left(
r^{\prime }\right) }\right) =0,
\end{equation*}%
there exists a process $Y=\left\{ Y_{t}:t\in \left[ 0,1\right] \right\} $
such that
\begin{equation}
Y_{0}=0\text{, \ }V\left( Y\right) <+\infty \text{ \ and \ }\lim_{r\uparrow
+\infty }V\left( Z^{\left( r\right) }-Y\right) =0,  \label{limvariation}
\end{equation}%
and moreover there exist a unique collection of kernels $f_{l,q}\in
L_{s}^{2}\left( \left[ 0,1\right] ^{l}\right) $ such that, for every $t\in %
\left[ 0,1\right] $, $Y_{t}$ admits the representation%
\begin{equation}
Y_{t}=\sum_{l=1}^{m+n}\sum_{\text{ }\left( l-n\right) \vee 1\leq q\leq
l\wedge m}I_{l}^{X}\left( \mathbf{1}_{A_{l,q}\left( t\right) }f_{l,q}\right)
,\text{ \ \ }t\in \left[ 0,1\right] \text{,}  \label{representation}
\end{equation}%
where, for every $k\geq 1$, we adopt the notation $\sum_{\text{
}k\leq q\leq 0}:=0.$ In particular, $Y\in \mathbf{FS}_{n+m-1}$.
\end{lemma}

\begin{proof}
If $m$ or $n$ is equal to zero, the statement can be proved by standard
arguments. Now suppose $n,m\geq 1$, and fix $r\geq 1$ and $u=1,...,M_{r}.$
The multiplication formula for multiple Wiener integrals yields%
\begin{equation*}
Z^{\left( u,r\right) }=\sum_{l=1}^{m+n}\sum_{\text{ }\left( l-n\right) \vee
1\leq q\leq l\wedge m}I_{l}^{X}\left( \widetilde{\left( h_{q}^{\left(
u,r\right) }\mathbf{1}_{\left( 0,t\right) }^{\otimes q}\right) \otimes
_{0}\left( g_{l-q}^{\left( u,r\right) }\mathbf{1}_{\left( t,1\right)
}^{\otimes l-q}\right) }\right)
\end{equation*}%
where $g_{0}^{\left( u,r\right) }:=b^{\left( u,r\right) }$ and~~$\widetilde{}
$~stands for symmetrization. Note that if $q=l$, then $l\leq m$ and%
\begin{equation*}
I_{l}^{X}\left( \widetilde{\left( h_{q}^{\left( u,r\right) }\mathbf{1}%
_{\left( 0,t\right) }^{\otimes q}\right) \otimes _{0}\left( g_{l-q}^{\left(
u,r\right) }\mathbf{1}_{\left( t,1\right) }^{\otimes l-q}\right) }\right)
=b^{\left( u,r\right) }I_{l}^{X}\left( \left( h_{l}^{\left( u,r\right) }%
\mathbf{1}_{A_{l,l}\left( t\right) }\right) \right) .
\end{equation*}%
On the other hand, when $1\leq q<l$, for every $\mathbf{x}_{l}\in \left[ 0,1%
\right] ^{l}$%
\begin{eqnarray*}
\widetilde{\left( h_{q}^{\left( u,r\right) }\mathbf{1}_{\left( 0,t\right)
}^{\otimes q}\right) \otimes _{0}\left( g_{l-q}^{\left( u,r\right) }\mathbf{1%
}_{\left( t,1\right) }^{\otimes l-q}\right) } &=&\binom{l}{q}^{-1}\sum_{%
\mathbf{j}_{\left( q\right) }\subset \{ 1,...,l\} }h_{q}^{\left(
u,r\right) }\left( \mathbf{x}_{\mathbf{j}_{\left( q\right)
}}\right) g_{l-q}^{\left( u,r\right) }\left( \mathbf{x}_{\{
1,...,l\}
\backslash \mathbf{j}_{\left( q\right) }}\right) \times \\
&&\text{ \ \ \ \ \ \ \ \ \ \ \ \ \ \ \ \ \ \ \ \ \ }\times \mathbf{1}_{\left[
0,t\right) ^{q}}\left( \mathbf{x}_{\mathbf{j}_{\left( q\right) }}\right)
\mathbf{1}_{\left( t,1\right] ^{l-q}}\left( \mathbf{x}_{\{ 1,...,l\}
\backslash \mathbf{j}_{\left( q\right) }}\right) \\
&=&\binom{l}{q}^{-1}\mathbf{1}_{A_{l,q}\left( t\right) }\left( \mathbf{x}%
_{l}\right) \times \\
&&\text{ \ \ \ \ \ \ \ \ \ }\times \sum_{\mathbf{j}_{\left( q\right)
}\subset \{ 1,...,l\}}h_{q}^{\left( u,r\right) }\left( \mathbf{x}_{%
\mathbf{j}_{\left( q\right) }}\right) g_{l-q}^{\left( u,r\right)
}\left( \mathbf{x}_{\{ 1,...,l\} \backslash \mathbf{j}_{\left(
q\right) }}\right) \mathbf{1}_{\Delta _{l}^{\mathbf{j}_{\left(
q\right) }}}\left( \mathbf{x}_{l}\right) .
\end{eqnarray*}

Since the function
\begin{equation*}
\mathbf{x}_{l}\mapsto \sum_{\mathbf{j}_{\left( q\right) }\subset
\{ 1,...,l\} }h_{q}^{\left( u,r\right) }\left(
\mathbf{x}_{\mathbf{j}_{\left( q\right)
}}\right) g_{l-q}^{\left( u,r\right) }\left( \mathbf{x}_{\{ 1,...,l\} \backslash \mathbf{%
j}_{\left( q\right) }}\right) \mathbf{1}_{\Delta
_{l}^{\mathbf{j}_{\left( q\right) }}}\left( \mathbf{x}_{l}\right)
\end{equation*}%
is symmetric, we immediately deduce that, for every $r\geq 1$, the family of
random variables
\begin{equation*}
\left\{ Z_{t}^{\left( r\right) }:t\in \left( 0,1\right) \right\} ,
\end{equation*}%
as defined in (\ref{approaching}), admits a representation of the form (\ref%
{representation}), and namely%
\begin{equation}
Z_{t}^{\left( r\right) }=\sum_{l=1}^{m+n}\sum_{\text{ }\left( l-n\right)
\vee 1\leq q\leq l\wedge m}I_{l}^{X}\left( \mathbf{1}_{A_{l,q}\left(
t\right) }f_{l,q}^{\left( r\right) }\right) ,  \label{expansion}
\end{equation}%
where
\begin{equation*}
f_{l,q}^{\left( r\right) }\left( \mathbf{x}_{l}\right) :=\binom{l}{q}%
^{-1}\sum_{u=1}^{M_{r}}\sum_{\mathbf{j}_{\left( q\right) }\subset \{ 1,...,l\}
}h_{q}^{\left( u,r\right) }\left( \mathbf{x}_{\mathbf{j}%
_{\left( q\right) }}\right) g_{l-q}^{\left( u,r\right) }\left( \mathbf{x}%
_{\{ 1,...,l\} \backslash \mathbf{j}_{\left( q\right) }}\right)
\mathbf{1}_{\Delta _{l}^{\mathbf{j}_{\left( q\right) }}}\left( \mathbf{x}%
_{l}\right) .
\end{equation*}%
Point (i) in the statement now follows from Theorem 6 and formula (\ref%
{expansion}). Now suppose that
\begin{equation*}
\lim_{r,r^{\prime }\rightarrow +\infty }V\left( Z^{\left( r\right)
}-Z^{\left( r^{\prime }\right) }\right) =0.
\end{equation*}%
Then, the existence of a process $Y$ satisfying (\ref{limvariation}) follows
from the same arguments contained in the proof of Theorem 4. Moreover,
relation (\ref{majoration}) implies immediately that for every $l$ and $q$,
the family $\left\{ f_{l,q}^{\left( r\right) }:r\geq 1\right\} $ is a Cauchy
sequence in $L_{s}^{2}\left( \left[ 0,1\right] ^{l}\right) $. Since $%
Y_{t}=L^{2}$-$\lim_{r\rightarrow +\infty }Z_{t}^{\left( r\right) }$ for
every $t$, the conclusion is obtained by standard arguments.
\end{proof}

\bigskip

Now, for every $p\geq 0$, call $\mathbf{BF}_{p}$ the subset of the class $%
\mathbf{BF}$, as defined through formula (\ref{FB}), composed of processes
with the form (\ref{approaching}) and such that $n+m\leq p$. We have
therefore the following

\bigskip

\begin{proposition}
Fix $N\geq 0$, and consider a measurable process $Y=\left\{ Y_{t}:t\in \left[
0,1\right] \right\} $. Then, the following conditions are equivalent:

\begin{enumerate}
\item $Y\in \mathbf{FS}_{N}$;

\item there exists a sequence $Z^{\left( r\right) }\in \mathbf{BF}_{N+1}$, $%
r\geq 1$, such that $\lim_{r\rightarrow +\infty }V\left( Z^{\left( r\right)
}-Y\right) =0$
\end{enumerate}
\end{proposition}

\begin{proof}
The implication 2. $\Longrightarrow $ 1. is an immediate consequence of
Lemma 7 and Theorem 6. To deal with the opposite direction, suppose that $%
Y_{t}=\delta \left( u\mathbf{1}_{\left[ 0,t\right] }\right) $, $t\in \left[
0,1\right] $, where $u_{\alpha }\left( \omega \right) \in L^{2}\left( \left[
0,1\right] \times \Omega \right) $ is such that, for every $\alpha \in \left[
0,1\right] ,$ $u_{\alpha }\in \oplus _{j=0,...,N}C_{j}$. Note that $u\in
\mathbb{L}^{k,p}$ for every $k\geq 1$ and $p>2$, and we can therefore take
up the same line of reasoning and notation as in the proof of Theorem 1. In
particular, according to Proposition 1 in Tudor (2004), we know that $%
Y $ admits the representation $Y_{t}=\int_{0}^{t}\mathbb{E}\left[ v_{\alpha
}\mid \mathcal{F}_{\left[ \alpha ,t\right] ^{c}}\right] dX_{\alpha }$, where
the process $v_{\alpha }=u_{\alpha }+\int_{0}^{\alpha }D_{\alpha
}u_{s}dX_{s} $, $\alpha \in \left[ 0,1\right] $, is also such that $%
v_{\alpha }\in \oplus _{j=0,...,N}C_{j}$ for every $\alpha .$ By linearity,
this implies that for every partition $\pi =\left\{
0=t_{0}<...<t_{n}=1\right\} $ the random variables $F_{i}^{\pi }$, $%
i=0,...,n-1$, as defined in (\ref{RVapproaching}), are such that $F_{i}^{\pi
}\in \oplus _{j=0,...,N}C_{j}$. According to the remark following Lemma 2,
every $F_{i}^{\pi }$ is the limit, say in $\mathbb{D}^{3,3}$, of a sequence
of random variables with the form%
\begin{equation*}
G_{m}^{\left( i,m\right) }=\sum_{k=1}^{M_{m}}G_{m,k}^{\left( i,\pi ,1\right)
}\times G_{m,k}^{\left( i,\pi ,2\right) }\text{, \ }m\geq 1\text{,}
\end{equation*}%
where $M_{m}\geq 1$ for every $m$, and also
\begin{eqnarray*}
G_{m,k}^{\left( i,\pi ,1\right) } &=&a+\sum_{l=1}^{\gamma
_{1}}I_{l}^{X}\left( h_{l}\mathbf{1}_{\left( 0,t_{j}\right) ^{l}}\right) \\
G_{m,k}^{\left( i,\pi ,2\right) } &=&b+\sum_{r=1}^{\gamma
_{2}}I_{r}^{X}\left( g_{r}\mathbf{1}_{\left( t_{j+1},1\right) ^{r}}\right)
\end{eqnarray*}%
where all dependencies on $i,\pi ,m$ and $k$ have been dropped in
the
second members, and $\gamma _{1}+\gamma _{2}\leq N+1$. By using relations (\ref%
{Bmart}) and (\ref{Bmart2}), we see immediately that the process $%
U_{t}^{\left( m,k,\pi ,i\right) }$, $t\in \left[ 0,1\right] $, is an element
of $\mathbf{BF}_{N+1}$, and the conclusion is obtained as in the proof of
Theorem 1.
\end{proof}

\section{Skorohod integrals as time-reversed Brownian martingales}

Now fix $k\geq3$ and $p>2$, take $u\in\mathbb{L}^{k,p}$, and note
$Y_{t}=\delta(u\mathbf{1}_{[0,t]})$. Suppose moreover
that the process $v_{\alpha }\in\mathbb{L}^{k-2,p}$ appearing in formula (\ref%
{TudorRepres}) is such that $v_{\alpha }=D_{\alpha }F$ for some
$F\in \mathbb{D}^{1,2}$ (we refer to Nualart (1995, p. 40) for a
characterization of such processes in term of their Wiener-It\^{o}
expansion). Then, according to the \textit{generalized Clark-Ocone
formula} stated in Nualart and Pardoux (1988),
\begin{equation}
Y_{t} =\int_{0}^{t}\mathbb{E}%
\left[ D_{\alpha }F\mid \mathcal{F}_{\left[ \alpha ,t\right]
^{c}}\right]
dX_{\alpha }=F-\mathbb{E}\left[ F\mid \mathcal{F}_{t^{c}}\right] ,\ \ t\in %
\left[ 0,1\right].  \label{NPa}
\end{equation}%
As made clear by the following discussion, a process of the type $Y_{t}=F-%
\mathbb{E}\left[ F\mid \mathcal{F}_{t^{c}}\right] $ can be easily
represented as a \textit{time-reversed Brownian martingale}. The principal
aim of this section is to establish sufficient conditions to have that $%
Y_{t} $ is a semimartingale in its own filtration (the reader is
referred to Tudor (2004), for further applications of (\ref{NPa})
to Skorohod integration).

\bigskip

To this end, for every $f\in L^{2}\left( \left[ 0,1\right] \right) $ we
define $\widehat{f}\left( x\right) =f\left( 1-x\right) $, so that the
transformation $f\mapsto \widehat{f}$ is an isomorphism of \thinspace $%
L^{2}\left( \left[ 0,1\right] \right) $ into itself. Such an
operator can be extended to the space $L_{s}^{2}\left( \left[
0,1\right] ^{n}\right) $
-- i.e. the space of square integrable and symmetric functions on $\left[ 0,1%
\right] ^{n}$ -- by setting%
\begin{equation*}
\widehat{f}_{n}\left( x_{1},...x_{n}\right) =f\left(
1-x_{1},...,1-x_{n}\right)
\end{equation*}%
for every $f_{n}\in L_{s}^{2}\left( \left[ 0,1\right] ^{n}\right) $, thus
obtaining an isomorphism of $L_{s}^{2}\left( \left[ 0,1\right] ^{n}\right) $
into itself. We also set, for $f\in L^{2}\left( \left[ 0,1\right] \right) $,
$\widehat{X}\left( f\right) =X\left( \widehat{f}\right) $ and eventually
\begin{equation*}
\widehat{X}=\left\{ \widehat{X}\left( f\right) :f\in L^{2}\left( \left[ 0,1%
\right] \right) \right\} .
\end{equation*}

Of course, $\widehat{X}$ is an isonormal Gaussian process on $L^{2}\left( %
\left[ 0,1\right] \right) $, and the random function%
\begin{equation*}
\widehat{X}_{t}=\widehat{X}\left( \mathbf{1}_{\left[ 0,t\right] }\right)
=X_{1}-X_{1-t}\text{, \ \ }t\in \left[ 0,1\right] ,
\end{equation*}%
is again a standard Brownian motion. As usual, given $n\geq 1$ and
$h_{n}\in L_{s}^{2}\left( \left[ 0,1\right] ^{n}\right) $,
$I_{n}^{X}\left( h_{n}\right) $ and $I_{n}^{\widehat{X}}\left(
h_{n}\right) $ stand for the
multiple Wiener-It\^{o} integrals of $h_{n}$, respectively with respect to $%
X $ and $\widehat{X}$ (see Nualart (1995)). The following lemma
will be useful throughout the sequel.

\begin{lemma}
Let $F\in L^{2}\left( \mathbb{P}\right) $ have the Wiener-It\^{o} expansion $%
F=\mathbb{E}\left( F\right) +\sum_{n=1}^{\infty }I_{n}^{X}\left(
f_{n}\right) $, then%
\begin{equation*}
F=\mathbb{E}\left( F\right) +\sum_{n=1}^{\infty }I_{n}^{\widehat{X}}\left(
\widehat{f}_{n}\right) .
\end{equation*}
\end{lemma}

\begin{proof}
By density, one can consider functionals with the form $F=I_{n}^{X}\left(
f^{\otimes n}\right) $, $n\geq 1$, where $f\in L^{2}\left( \left[ 0,1\right]
\right) $ and $f^{\otimes n}\left( x_{1},...x_{n}\right) =f\left(
x_{1}\right) ...f\left( x_{n}\right) $. In this case, it is well known that $%
F=n!H_{n}\left( X\left( f\right) \right) $, where $H_{n}$ is the
$n$th Hermite polynomial as defined in Nualart (1990, Ch. 1), and
therefore
\begin{equation*}
F=n!H_{n}\left( \widehat{X}\left( \widehat{f}\right) \right) =I_{n}^{%
\widehat{X}}\left( \widehat{f}^{\otimes n}\right) =I_{n}^{\widehat{X}}\left(
\widehat{f^{\otimes n}}\right) ,
\end{equation*}%
thus proving the claim.
\end{proof}

\bigskip

We now introduce the following filtration:
\begin{equation*}
\widehat{\mathcal{F}}_{t}=\sigma \left\{ \widehat{X}_{h}:h\leq t\right\}
\text{, \ \ }t\in \left[ 0,1\right] .
\end{equation*}

Note that%
\begin{eqnarray}
\mathcal{F}_{\left[ s,t\right] ^{c}} &=&\mathcal{F}_{s}\vee \widehat{%
\mathcal{F}}_{1-t}  \label{sigmafields} \\
\mathcal{F}_{t^{c}} &=&\widehat{\mathcal{F}}_{1-t}.  \notag
\end{eqnarray}

\begin{proposition}
Let $\left\{ Y_{t}:t\in \left[ 0,1\right] \right\} $ be a measurable process.

\begin{enumerate}
\item The following conditions are equivalent,
\end{enumerate}

\begin{description}
\item[(i)] there exists $F\in L^{2}\left( \mathbb{P}\right) $ such that $%
Y_{t}=F-\mathbb{E}\left( F\mid \mathcal{F}_{t^{c}}\right) ;$

\item[(ii)] there exists a square integrable $\widehat{\mathcal{F}}_{t}$ - martingale $%
\left\{ \widehat{M}_{t}:t\in \left[ 0,1\right] \right\} $ such that $Y_{t}=%
\widehat{M}_{1}-\widehat{M}_{1-t}$;

\item[(iii)] there exists a $\widehat{\mathcal{F}}_{\alpha }$ - predictable
process $\left\{ \widehat{\phi }_{\alpha }:\alpha \in \left[ 0,1\right]
\right\} $, such that $\mathbb{E}\left( \int_{0}^{1}\widehat{\phi }_{\alpha
}^{2}d\alpha \right) <+\infty $ and \ $Y_{t}=\int_{1-t}^{1}\widehat{\phi }%
_{\alpha }d\widehat{X}_{\alpha };$

\item[(iv)] there exist kernels $f_{n}\in L_{s}^{2}\left( \left[ 0,1\right]
^{n}\right) $, $n\geq 1$, such that
\begin{equation*}
Y_{t}=\sum_{n=1}^{\infty }I_{n}^{X}\left( f_{n}\left( 1-\mathbf{1}_{\left[
t,1\right] }^{\otimes n}\right) \right) =\sum_{n=1}^{\infty }I_{n}^{\widehat{%
X}}\left( \widehat{f}_{n}\left( 1-\mathbf{1}_{\left[ 0,1-t\right]
}^{\otimes n}\right) \right) ,
\end{equation*}
where the convergence of the series takes place in $ L^{2}\left(
\mathbb{P}\right) $.
\end{description}

\begin{enumerate}
\item[2.] Let either one of conditions (i)-(iv) be verified, and let $F$ be
given by (i) and the $f_{n}$'s by (iv). Then,%
\begin{equation*}
F=\mathbb{E}\left( F\right) +\sum_{n=1}^{\infty }I_{n}^{X}\left(
f_{n}\right) =\mathbb{E}\left( F\right) +\sum_{n=1}^{\infty }I_{n}^{\widehat{%
X}}\left( \widehat{f}_{n}\right) .
\end{equation*}

\item[3.] Under the assumptions of point 2, suppose moreover that
$F$ is an element of $\mathbb{D}^{1,2}$, and let $\widehat{\phi }$
be given by
(iii). Then,%
\begin{equation}
\widehat{\phi }_{\alpha }=\mathbb{E}\left[ D_{1-\alpha }F\left( X\right)
\mid \widehat{\mathcal{F}}_{\alpha }\right] ,\text{ \ }\alpha \in \left[ 0,1%
\right] \text{,}\label{wu}
\end{equation}%
where $DF\left( X\right) $ is the usual Malliavin derivative of $F$,
regarded as a functional of $X$.
\end{enumerate}
\end{proposition}

\textbf{Remark} -- Note that formula (\ref{wu}) above appears also
in Wu (1990, formula (4.4)), where it is obtained by completely
different arguments.
\medskip

\begin{proof}
If (i) is verified, then (ii) holds, thanks to (\ref{sigmafields}), by
defining $\widehat{M}_{t}=\mathbb{E}\left( F\mid \widehat{\mathcal{F}}%
_{t}\right) $. On the other hand, (ii) implies (iii) due to the predictable
representation property of $\widehat{X}$. Of course, if (iii) is verified,
then%
\begin{equation*}
Y_{t}=\int_{1-t}^{1}\widehat{\phi }_{\alpha }d\widehat{X}_{\alpha
}=\int_{0}^{1}\widehat{\phi }_{\alpha }d\widehat{X}_{\alpha }-\int_{0}^{1-t}%
\widehat{\phi }_{\alpha }d\widehat{X}_{\alpha }=F-\mathbb{E}\left[ F\mid
\widehat{\mathcal{F}}_{1-t}\right] ,
\end{equation*}%
where $F=\int_{0}^{1}\widehat{\phi }_{\alpha }d\widehat{X}_{\alpha }$, thus
proving the implication (iii) $\Longrightarrow $ (i). Now, let (i) be
verified, and let $F$ have the representation
\begin{equation*}
F=\mathbb{E}\left( F\right) +\sum_{n=1}^{\infty }I_{n}^{X}\left(
f_{n}\right) ;
\end{equation*}%
we may apply Lemma 1.2.4 in Nualart (1995) to obtain that%
\begin{equation}
Y_{t}=\sum_{n=1}^{\infty }I_{n}^{X}\left( f_{n}\right) -\mathbb{E}\left[
\sum_{n=1}^{\infty }I_{n}^{X}\left( f_{n}\right) \mid \mathcal{F}_{t^{c}}%
\right] =\sum_{n=1}^{\infty }I_{n}^{X}\left( f_{n}\right)
-\sum_{n=1}^{\infty }I_{n}^{X}\left( f_{n}\mathbf{1}_{\left[ t,1\right]
}^{\otimes n}\right)  \label{conditioning}
\end{equation}%
thus giving immediately (i) $\Longrightarrow $ (iv) (the second equality in
(iv) is a consequence of Lemma 9). The opposite implication may be obtained
by reading backwards formula (\ref{conditioning}). The proof of point 2 is
now immediate. To deal with point 3, observe if $F$ is derivable in the
Malliavin sense as a functional of $X$, then $F$ is also derivable as a
functional of $\widehat{X}$, and the two derivative processes must verify%
\begin{equation*}
D_{\alpha }F\left( \widehat{X}\right) =D_{1-\alpha }F\left( X\right) ,\text{
\ \ a.e. -- }d\alpha \otimes d\mathbb{P}\text{,}
\end{equation*}%
where $DF\left( \widehat{X}\right) $ stands for the Malliavin derivative of $%
F$, regarded as a functional of $\widehat{X}$. As a matter of fact, let $%
F_{k}$ be a sequence of polynomial functionals with the form $F_{k}=p\left(
X\left( h_{1}\right) ,...,X\left( h_{m}\right) \right) $, where $p$ is a
polynomial in $m$ variables (note that $p$, $m$ and the $h_{j}$'s may in
general depend on $k$), converging to $F$ in $L^{2}\left( \mathbb{P}\right) $
and satisfying%
\begin{equation*}
\mathbb{E}\left[ \int_{0}^{1}\left( \sum_{j=1}^{m}\frac{\partial }{\partial
x_{j}}p\left( X\left( h_{1}\right) ,...,X\left( h_{m}\right) \right)
h_{j}\left( x\right) -D_{x}F\left( X\right) \right) ^{2}dx\right]
\rightarrow 0\text{. }
\end{equation*}

Then, $p\left( X\left( h_{1}\right) ,...,X\left( h_{m}\right) \right)
=p\left( \widehat{X}\left( \widehat{h}_{1}\right) ,...,\widehat{X}\left(
\widehat{h}_{m}\right) \right) $, and also
\begin{equation*}
\mathbb{E}\left[ \int_{0}^{1}\left( \sum_{j=1}^{m}\frac{\partial }{\partial
x_{j}}p\left( \widehat{X}\left( \widehat{h}_{1}\right) ,...,\widehat{X}%
\left( \widehat{h}_{m}\right) \right) \widehat{h}_{j}\left( x\right)
-D_{1-x}F\left( X\right) \right) ^{2}dx\right] \rightarrow 0\text{,}
\end{equation*}%
thus giving immediately the desired conclusion. The proof of point
3 is achieved by using the Clark-Ocone formula (see Clark (1970)
and Ocone (1984)).
\end{proof}

\bigskip

\textbf{Example --} Let $F=H_{n}\left( X\left( h\right) \right) $, where $%
H_{n}$ is the $n$th Hermite polynomial and $h$ is such that $\left\Vert
h\right\Vert =1$. Then, thanks to Proposition 10-2 the process $Y_{t}=F-%
\mathbb{E}\left[ F\mid \mathcal{F}_{t^{c}}\right] $ has the representation%
\begin{eqnarray}
Y_{t} &=&\frac{1}{n!}\left[ I_{n}^{X}\left( h^{\otimes n}\right)
-I_{n}^{X}\left( h^{\otimes n}\mathbf{1}_{\left[ t,1\right] }^{\otimes
n}\right) \right] =H_{n}\left( X\left( h\right) \right) -\left\Vert h\mathbf{%
1}_{\left[ 0,t\right] }\right\Vert ^{n}I_{n}^{X}\left( \left( \frac{h\mathbf{%
1}_{\left[ t,1\right] }}{\left\Vert h\mathbf{1}_{\left[ t,1\right]
}\right\Vert }\right) ^{\otimes n}\right)  \label{hermite} \\
&=&H_{n}\left( X\left( h\right) \right) -\left\Vert h\mathbf{1}_{\left[ t,1%
\right] }\right\Vert ^{n}H_{n}\left( \frac{X\left( h\mathbf{1}_{\left[ t,1%
\right] }\right) }{\left\Vert h\mathbf{1}_{\left[ t,1\right] }\right\Vert }%
\right)  \notag
\end{eqnarray}%
as well as
\begin{equation*}
Y_{t}=\int_{1-t}^{1}\mathbb{E}\left[ H_{n-1}\left( X\left( h\right) \right)
\mid \widehat{\mathcal{F}}_{\alpha }\right] h\left( 1-\alpha \right) d%
\widehat{X}_{\alpha }.
\end{equation*}

Formula (\ref{hermite}) generalizes the obvious relations (corresponding to
the case $n=1$ and $h=\mathbf{1}_{\left[ 0,1\right] }$)%
\begin{equation*}
X_{1}-\mathbb{E}\left[ X_{1}\mid \mathcal{F}_{t^{c}}\right] =X_{t}=\widehat{X%
}_{1}-\widehat{X}_{1-t}
\end{equation*}

\bigskip

Given a filtration $\left\{ \mathcal{G}_{t}:t\in \left[ 0,1\right] \right\} $%
, and two adapted, cadlag processes $U_{t}$, and $V_{t}$, we will write $%
\left[ U,V\right] =\left\{ \left[ U,V\right] _{t}:t\in \left[ 0,1\right]
\right\} $ to indicate the \textit{quadratic covariation process} of $U$ and
$V$ (if it exists). This means that $\left[ U,V\right] $ is the cadlag $%
\mathcal{G}_{t}$ - adapted process of bounded variation such that, for every
$t\in \left[ 0,1\right] $ and for every sequence of (possibly random)
partitions of $\left[ 0,t\right] $ -- say $\tau _{n}=\left\{
0<t_{1,n}<...<t_{M_{n},n}=t\right\} $ -- with mesh tending to zero, the
sequence%
\begin{equation*}
\lim_{n}\left[ U_{0}V_{0}+\sum_{i=0}^{M_{n}-1}\left(
U_{t_{i+1,n}}-U_{t_{i,n}}\right) \left( V_{t_{i+1,n}}-V_{t_{i,n}}\right) %
\right] =\left[ U,V\right] _{t}
\end{equation*}%
where the convergence is in probability, and uniform on compacts. The next
result uses quadratic covariations to characterize processes of the form $%
t\mapsto \left( F-\mathbb{E}\left[ F\mid \mathcal{F}_{t^{c}}\right] \right) $
in terms of semimartingales.

\bigskip

\begin{proposition}
Let $F$ and $\left\{ Y_{t}:t\in \left[ 0,1\right] \right\} $ satisfy either
one of conditions (i)-(iv) in Proposition 10, fix $k\geq 1$, and let $%
\widehat{\phi }_{\alpha }$, as in Proposition 10-1-(iii), be
c\`adl\`ag and of
the form%
\begin{equation*}
\widehat{\phi }_{\alpha }=\Phi \left( \alpha ;\widehat{X}\left( g_{1}\mathbf{%
1}_{\left[ 0,\alpha \right] }\right) ,...,\widehat{X}\left( g_{k}\mathbf{1}_{%
\left[ 0,\alpha \right] }\right) \right)
\end{equation*}%
where $\Phi $ is a measurable function on $\left[ 0,1\right] \times \Re ^{k}$%
, and $g_{j}\in L^{2}\left( \left[ 0,1\right] \right) $, $j=1,...,k$. If
there exists the quadratic covariation process $\left[ \widehat{\phi },%
\widehat{X}\right] $, then $Y_{t}$ is a semimartingale on $\left[ 0,1\right]
$ in its own filtration, and moreover%
\begin{equation}
Y_{t}=\int_{0}^{t}\widehat{\phi }_{1-\alpha }dX_{\alpha }-\left[ \widehat{%
\phi },\widehat{X}\right] _{1}+\left[ \widehat{\phi },\widehat{X}\right]
_{1-t}.  \label{egalite}
\end{equation}
\end{proposition}

\begin{proof}
The proof is directly inspired by Theorem 3.3\ in Jacod and
Protter (1988). Let $t\in \left( 0,1\right] $ and $\tau =\left\{
1-t=s_{0}<...<s_{n}=1\right\} $ be a deterministic partition of $\left[ 1-t,1%
\right] $. Then, when the mesh of $\tau $ converges to zero, $Y_{t}$ is
(uniformly) the limit in probability of
\begin{equation*}
\sum_{i=0}^{n-1}\widehat{\phi }_{s_{i}}\left( \widehat{X}_{s_{i+1}}-\widehat{%
X}_{s_{i}}\right) .
\end{equation*}

Now note that, since $\widehat{X}\left( g_{j}\mathbf{1}_{\left[ 0,1-\alpha %
\right] }\right) =X\left( \widehat{g}_{j}\right) -X\left( \widehat{g}_{j}%
\mathbf{1}_{\left[ 0,\alpha \right] }\right) $, $j=1,...,k$, the process $%
\alpha \mapsto \widehat{\phi }_{1-\alpha }$ is left-continuous and adapted
to the filtration%
\begin{equation*}
\mathcal{H}_{\alpha }=\sigma \left( X_{h}\text{, \ }h\leq \alpha \right)
\vee \sigma \left( X_{1},X\left( \widehat{g}_{1}\right) ,...,X\left(
\widehat{g}_{k}\right) \right) \text{, \ \ }\alpha \in \left[ 0,1\right] .
\end{equation*}

Therefore, since $X_{t}$ is classically a $\mathcal{H}_{t}$ - semimartingale
(see Chaleyat-Mauriel and Jeulin (1983)), the stochastic integral in (%
\ref{egalite}) is well defined as the limit in probability of the sequence%
\begin{equation*}
\sum_{i=0}^{n-1}\widehat{\phi }_{1-t_{i+1}}\left(
X_{t_{i}}-X_{t_{i+1}}\right) =\sum_{i=0}^{n-1}\widehat{\phi }%
_{s_{i+1}}\left( \widehat{X}_{s_{i+1}}-\widehat{X}_{s_{i}}\right)
\end{equation*}%
where $t_{i}=1-s_{i}$. Eventually, we shall observe that the finite
variation process $t\mapsto \left[ \widehat{\phi },\widehat{X}\right] _{1}-%
\left[ \widehat{\phi },\widehat{X}\right] _{1-t}$ is by definition the limit
in probability (as the mesh of $\tau $ converges to zero) of
\begin{equation*}
\sum_{i=0}^{n-1}\left( \widehat{\phi }_{s_{i+1}}-\widehat{\phi }%
_{s_{i}}\right) \left( \widehat{X}_{s_{i+1}}-\widehat{X}_{s_{i}}\right) ,
\end{equation*}%
and therefore it is a $\mathcal{H}_{t}$ - semimartingale, being an adapted
process of finite variation (to prove the adaptation, just observe that if $%
1-t\leq s\leq 1$, then
\begin{eqnarray*}
\widehat{\phi }_{s} &=&\Phi \left( \alpha ;X\left( \widehat{g}_{1}\right)
-X\left( \widehat{g}_{1}\mathbf{1}_{\left[ 0,1-s\right] }\right)
,...,X\left( \widehat{g}_{k}\right) -X\left( \widehat{g}_{k}\mathbf{1}_{%
\left[ 0,1-s\right] }\right) \right) \\
&\in &\sigma \left( X_{h}\text{, \ }h\leq t\right) \vee \sigma \left(
X_{1},X\left( \widehat{g}_{1}\right) ,...,X\left( \widehat{g}_{k}\right)
\right) ).
\end{eqnarray*}%
As a consequence of the above discussion, the quantity
\begin{equation*}
Y_{t}-\int_{0}^{t}\widehat{\phi }_{1-\alpha }dX_{\alpha }+\left[ \phi ,%
\widehat{X}\right] _{1}-\left[ \phi ,\widehat{X}\right] _{1-t}
\end{equation*}%
is the limit in probability of
\begin{equation*}
\sum_{i=0}^{n-1}\widehat{\phi }_{s_{i}}\left( \widehat{X}_{s_{i+1}}-\widehat{%
X}_{s_{i}}\right) -\sum_{i=0}^{n-1}\widehat{\phi }_{s_{i+1}}\left( \widehat{X%
}_{s_{i+1}}-\widehat{X}_{s_{i}}\right) +\sum_{i=0}^{n-1}\left(
\widehat{\phi
}_{s_{i+1}}-\widehat{\phi }_{s_{i}}\right) \left( \widehat{X}_{s_{i+1}}-%
\widehat{X}_{s_{i}}\right)
\end{equation*}%
which equals zero for every $\tau $. To conclude, observe that $Y_{t}$ is
the sum of two $\mathcal{H}_{t}$ - semimartingales, and therefore it is
itself a $\mathcal{H}_{t}$ - semimartingale and consequently, by Stricker's
theorem, it is a semimartingale in its own filtration.
\end{proof}

\bigskip

Now we state a (classic) sufficient condition for the existence of the
quadratic covariation process $\left[ \widehat{\phi },\widehat{X}\right] .$

\bigskip

\begin{proposition}
Under the assumptions and notation of Proposition 11, suppose that the
function $\Phi $ is of class $C^{1}$ in $\left[ 0,1\right] \times \Re ^{k}$.
Then, the quadratic covariation process $\left[ \widehat{\phi },\widehat{X}%
\right] $ exists.
\end{proposition}

\begin{proof}
This is an application of Theorem 5 in Meyer (1976, p. 359). The
vector
\begin{equation*}
\mathbf{\gamma }_{\alpha }:=\left( \alpha ,\widehat{X}_\alpha,\widehat{X}\left( g_{1}\mathbf{1}%
_{\left[ 0,\alpha \right] }\right) ,...,\widehat{X}\left( g_{k}\mathbf{1}_{%
\left[ 0,\alpha \right] }\right) \right)
\end{equation*}%
is indeed a $\left( k+2\right) $ - dimensional $\widehat{\mathcal{F}}%
_{\alpha }$ - semimartingale. Now define
\[ \Phi ^{\ast }\left(
\alpha,x_{1},...,x_{k+1}\right) =\Phi \left(
\alpha,x_{2},...,x_{k+1}\right) \text{, \ \ }\left(
\alpha,x_{1},...,x_{k+1}\right) \in \left[ 0,1\right] \times \Re
^{k+1}.
\]
Since the assumptions imply that $\Phi ^{\ast }$ is of class
$C^{1}$ in $\left[ 0,1\right] \times \Re ^{k+1}$
and $\widehat{\phi }%
_{\alpha }=\Phi ^{\ast } \left( \mathbf{\gamma }_{\alpha }\right)
$, the quadratic variation process $\alpha \mapsto \left[ \widehat{\phi },%
\widehat{\phi }\right] _{\alpha }$ exists, as well as the
processes $\left[ \widehat{X},\widehat{X}\right] $ and $\left[
\widehat{\phi }+\widehat{X},\widehat{\phi }+\widehat{X}\right] $. It follows that $\left[ \widehat{\phi },%
\widehat{X}\right] $ exists, thanks to the polarization identity%
\begin{equation*}
\left[ \widehat{\phi },\widehat{X}\right] _{\alpha }=\frac{1}{2}\left\{ %
\left[ \widehat{\phi }+\widehat{X},\widehat{\phi }+\widehat{X}\right]
_{\alpha }-\left[ \widehat{X},\widehat{X}\right] _{\alpha }-\left[ \widehat{%
\phi },\widehat{\phi }\right] _{\alpha }\right\} \text{, \ \ }\alpha \in %
\left[ 0,1\right] \text{.}
\end{equation*}
\end{proof}

\section{Anticipating integrals and stopping times}

For the sake of completeness, in this section we explore some
links between Skorohod integral processes and the family of
stopping times. Classically, the stopping times are strongly
related to the martingale theory. For instance, fix a filtration
$\mathcal{U}_{t}$ as well as a $\mathcal{U}_{t}$ - stopping time
$T$: it is well known, from the \emph{Optional Sampling Theorem
}(see e.g. Chung (1974)),\emph{\ }that, for any $\mathcal{U}_{t}$
- martingale $M_{t}$, the stopped process $t\mapsto M_{T\wedge t}$
is again a martingale for the filtration $t\mapsto
{\mathcal{U}}_{T\wedge t}$ of events determined prior to $T$. It
is also well-known that, a stopped It\^{o} integral at the
stopping time $T$ coincides with the It\^{o} integral on the
random interval $[0,T]$. In this section, we prove a variant of
the Optional Sampling Theorem for Skorohod integral processes and
we discuss what happens if one samples such a\ process at a random
time. For a discussion in this direction, see also the paper
Nualart and Thieullen (1994). We keep the notation of the previous
sections, and consider anticipating integral processes given by
\begin{equation*}
Y_{t}=\delta \left( u\mathbf{1}_{[0,t]}(\cdot )\right)
\end{equation*}%
where $u\mathbf{1}_{\left[ 0,t\right] }$ belongs to $Dom(\delta )$ for every
$t\in \left[ 0,1\right] $. Given two stopping times $S,$ $T$ for the
filtration ${\mathcal{F}}_{t}$, we denote by ${\mathcal{F}}_{T}$, resp. ${%
\mathcal{F}}_{S}$, the $\sigma $-field of the events determined prior to $T$%
, resp. $S$.

\smallskip

We have the following \emph{Optional Sampling Theorem.}

\begin{proposition}
If $S,T$ are ${\mathcal{F}}_{t}$- stopping times such that $S\leq T$ a.s.,
it holds that
\begin{equation}  \label{eq1}
E\left[ Y_{T}-Y_{S}|{\mathcal{F}}_{S}\right] =0.
\end{equation}
\end{proposition}

\begin{proof}
Let us first consider as in Karatzas and Shreve (1991) two sequences of stopping times $%
(S_{n})_{n},(T_{n})_{n}$ taking on a countable number of values in the
dyadic partition of $[0,1]$ and such that $S_{n}\rightarrow S$, $%
T_{n}\rightarrow T$ and
\begin{equation*}
S\leq S_{n}\mbox{  ,  }T\leq T_{n}\mbox{ and }S_{n}\leq T_{n}.
\end{equation*}%
As in Chung (1974), p. 325, using the fact that the process
$\left( \mathbb{E}\left( Y_{t}|{\mathcal{F}}_{t}\right) \right)
_{t}$ is a
martingale, we can prove that $\int_{A}Y_{S_{n}}d\mathbb{P}%
=\int_{A}Y_{T_{n}}d\mathbb{P}$ for every $A\in
\mathcal{F}_{S_{n}}$. We follow next the lines of the proof of
Theorem 1.3.22 in Karatzas and Shreve (1991), observing that the
sequence $(Y_{S_{n}})_{n}$ is uniformly integrable. This is
consequence of the bound
\begin{equation*}
\sup_{t}\mathbb{E}Y_{t}^{2}\leq \sup_{t}\left( \mathbb{E}(Y_{1}-Y_{t})^{2}+%
\mathbb{E}Y_{t}^{2}\right) \leq V(Y).
\end{equation*}
\end{proof}

\bigskip

The next result is a version of Theorem 2.5 of Nualart and
Thieullen (1994).

\begin{proposition}
Let $u\in \mathbb{L} ^{1,p}$, $p>4$, and let $T$ be a stopping
time for the
filtration ${\mathcal{F}}_{t}$. Then $u\mathbf{1}_{[0,T]}$ belongs to $%
Dom(\delta )$ and it holds
\begin{equation}  \label{stop}
\delta (u\mathbf{1}_{[0,t]} ) \mid _{t=T} = \delta (u\mathbf{1}_{[0,T]}).
\end{equation}
\end{proposition}

 \begin{proof} Since, for $u$ as in the statement, the process $t\mapsto \int_{0}^{t}\mathbb{%
E}\left( u_{s}\right) dX_{s}$ is a continuous, square integrable Gaussian $\mathcal{F}%
_{t}$ - martingale, we can assume, without loss of generality, that $\mathbb{E%
}\left( u_{t}\right) =0$ for every $t\in \left[ 0,1\right] $. We
first prove  the property (\ref{stop}) for the
 approximation $u^{\pi }$ given by (\ref{CondStrato})
 $$u^{\pi }_{t}= \sum _{i=0}^{n-1}\frac{1}{t_{i+1}-t_{i}}\left(\int _{t_{i}}^{t_{i+1}}E\left( u_{s}
 \mid {\cal{F}}_{[t_{i},t_{i+1}]^{c}}\right) ds\right)
 \mathbf{1}_{[t_{i},t_{i+1}]}(t).$$ Let us consider the sum
 $$S=\sum _{i=0}^{n-1} F_{i}\left( X_{T\wedge
 t_{i+1}}-X_{T\wedge t_{i}}\right) = \sum _{i=0}^{n-1} F_{i} \delta
 (\mathbf{1}_{[0,T]} \mathbf{1}_{[t_{i}, t_{i+1}]})$$
 where $F_{i}=\frac{1}{t_{i+1}-t_{i}}\left(\int _{t_{i}}^{t_{i+1}}E\left( u_{s}
 \mid {\cal{F}}_{[t_{i},t_{i+1}]^{c}}\right) ds\right)$.
 Using relation (\ref{ip}) (note that all hypothesis  are
 satisfied, that is, $F_{i}\in \mathbb{D}^{1,2}$, $\mathbf{1}_{[0,T]}
 \mathbf{1}_{[t_{i}, t_{i+1}]}\in Dom(\delta ),$ being adapted, and $\mathbb{E}\left( F^{2} \int
 _{0}^{1}\mathbf{1}_{[0,T]}(s)
 \mathbf{1}_{[t_{i}, t_{i+1}]}(s)ds\right) \leq \mathbb{E}(F^{2})<\infty $ ) and (\ref{d}),
  we obtain that $u^{\pi } \mathbf{1}_{[0,T]}\in Dom (\delta )$ and
 $$\delta (u^{\pi } \mathbf{1}_{[0,T]}) =S=\sum _{i=0}^{n-1} F_{i}\left( X_{t\wedge
 t_{i+1}}-X_{t\wedge t_{i}}\right)  \mid _{t=T} =  \delta (u^{\pi }\mathbf{1}_{[0,t]} ) \mid _{t=T} .$$

 Now recall that, for every partition $\pi $, the process $u^{\pi }$ is an
element of $\mathbb{L}^{1,p}$, and also, when $\left\vert \pi
\right\vert
\rightarrow 0$,%
\begin{equation}
\begin{array}{ll}
u^{\pi }\rightarrow u & \text{in }\mathbb{L}^{1,p} \\
u^{\pi }\mathbf{1}_{\left[ 0,T\right] }\rightarrow u\mathbf{1}_{\left[ 0,T%
\right] } & \text{in }L^{2}\left( \left[ 0,1\right] \times \Omega
\right)
\\
\delta \left( u^{\pi }\mathbf{1}_{\left[ 0,t\right] }\right)
\rightarrow
\delta \left( u\mathbf{1}_{\left[ 0,t\right] }\right)  & \text{in }%
L^{2}\left( \mathbb{P}\right) \text{ for every }t\in \left[ 0,1\right] .%
\end{array}
\label{Convrg1}
\end{equation}

Fix a sequence of partitions $\pi $ such that $\left\vert \pi
\right\vert \rightarrow 0$. From (\ref{Convrg1}), we deduce
immediately that there exists a finite constant $K>0$, not
depending on $\pi $, such that
\[
\int_{0}^{1}\mathbb{E}\left[ \left\vert \int_{0}^{1}\left(
D_{s}u_{t}^{\pi
}\right) ^{2}ds\right\vert ^{\frac{p}{2}}\right] dt<K\text{, \ \ for every }%
\pi \text{.}
\]

Moreover, since $\mathbb{E}\left( u_{t}^{\pi }\right) =0$ for
every $t$, we can use the same line of reasoning as in the proof
of Nualart (1998, Proposition 5.1.1), and deduce the existence of
a finite constant $K'>0$ such that, for every $s,t\in \left[
0,1\right] $ and every $\pi$,
\[
\mathbb{E}\left[ \left\vert \delta \left( u^{\pi }\mathbf{1}_{\left[ 0,t%
\right] }\right) -\delta \left( u^{\pi }\mathbf{1}_{\left[
0,s\right] }\right) \right\vert ^{p}\right] \leq K' \times
\left\vert t-s\right\vert ^{\frac{p}{2}-1}.
\]

As a consequence, by applying for instance Nualart (1998, Lemma
5.3.1), and since $T$ takes values in $\left[ 0,1\right] $ by
construction, we deduce that, as $\left\vert \pi \right\vert
\rightarrow 0$,
\[
\delta \left( u^{\pi }\mathbf{1}_{\left[ 0,T \right] }\right)
=\left. \delta \left( u^{\pi }\mathbf{1}_{\left[ 0,t\right]
}\right)
\right\vert _{t=T}\rightarrow \left. \delta \left( u\mathbf{1}_{\left[ 0,t%
\right] }\right) \right\vert _{t=T}\text{ \ \ in }L^{p}\left( \mathbb{P}%
\right) .
\]

We conclude by the basic lemma for the convergence of Skorohod
integrals that $u\mathbf{1}_{[0,T]}\in Dom (\delta )$ and
(\ref{stop}) holds.
\end{proof}

\smallskip

\textbf{Remark} -- Note  that in Nualart and Thieullen (1994,
Theorem 2.5) the authors proved the following relation, for every
${\cal{F}}_{t}$-stopping time $T$ and for every $u\in Dom (\delta
)$,
$$\delta (u\mathbf{1}_{[0,T]}) =\delta (u\mathbf{1}_{[0,t]} )
\mid _{t=T^{+}} $$ where $\delta (u\mathbf{1}_{[0,t]} ) \mid
_{t=T^{+}} $ is defined as
$$\delta (u\mathbf{1}_{[0,t]} ) \mid
_{t=T^{+}}=\lim _{\varepsilon \to 0} \frac{1}{\varepsilon } \int
_{T}^{T+\varepsilon } \delta (u\mathbf{1} _{[0,s]} ) ds $$ when
the above limit exists in $L^{2}(\mathbb{P})$. The obtention of
the result (\ref{stop}) is due to the use of the approximating
processes  (\ref{CondStrato}) for which the limit can be
explicitly computed. Note that, with our method, we do not need to
introduce any special assumption on $T$. On the other hand, we are
forced to assume a stronger hypothesis on the integrand $u$, that
is, $u\in \mathbb{L}^{1,p}$, $p>4$, instead of $u\in Dom (\delta
)$.

\end{document}